\documentclass[a4paper,10pt]{article}
\usepackage[utf8]{inputenc}
\usepackage{amsmath,amssymb}

\begin{document}
\title{Some remarks on monodromy}
\author{Tove Dahn}
\maketitle
\begin{abstract}
We consider hypoelliptic symbols, over a very regular Lie group and discuss monodromy for a spectral stratification, using results of Nilsson and B\"acklund.
\end{abstract}

\section{Conjugation}
Starting from $dU \in G$, define $G_{r}$ by dV, so that $<Uf,Vg>=<f,g>$ for $f \in (I)$, $g \in (J)$,
that is conjugation relative a bilinear form over ideals with topology for instance H.
We can define a zero set $\subset \cup_{j} E_{j}$ relative the bilinear form, such that $Vf \subset (I)(E_{j})$ and $\mid   {}^{t} UVf \mid < \epsilon$ over $(I)(E_{j})$ (normalized), assuming dU,dV BV measures (of bounded variation).
A zero function has support on a zero set. We can in this way define (f,g), a class of zero functions to $dU \times dV$.
ex: $g=\Phi(f)$ relative $<,>$, such that $U \overline{f}=V\Phi f$, thus conjugation defines a homogeneous space.

Consider $< Uf,Vg>=<f,g>$ as volume preserving (relative scaling) conjugation.
When g is dependent of f according to g=g(f) (homogeneity),
we have conjugation according to $Vg(f)=g(Uf)$.
Assume $\Phi_{U}(f)=\int f dU + \int_{\Gamma(U)} df$.
The dimension for $d V^{-1}U \simeq d U^{-1}V$, refers to the domain where $\Phi_{V/U}f$ is regular iff $\Phi_{U/V}f$ is regular.
$d V^{-1}U \rightarrow d U^{-1}V$ bi-continuous, defines an abelian group.

Consider $dV=\rho dU$ as continuation, when $\rho=const$ given change of type, we have presence of trace, for instance a spiral axis. When the type is preserved, we have a non-trivial polar, that is not dense range.

Consider the flux condition, $\Phi(U)(f)-\Phi(V)(f)=\int_{\Gamma(U)} df - \int_{\Gamma(V)} df =0$, corresponding to $\Gamma(U) \sim \Gamma(V)$. When the representation has no trace, then we must have
$d V^{-1}U=dI$ implies $d U=dI$. We consider conjugation relative rectifiable trajectories $\Gamma$, with $\{ dU=dI \} \rightarrow \Gamma$ absolute continuous, . When dU and the harmonic conjugate $d U^{\diamondsuit}$ are both closed, harmonic conjugation maps on to complex conjugation. Given $dU \in G$ and  $<Uf,Vg>=<f,g>$, as ${}^{t}VU=I$ and $dI \in G$, we have $d {}^{t}  V \in G$.
When U idempotent ${}^{t} U U =I$, that is given maximal dimension for symmetric regularity, we have $V^{-1}UV=U$. An elastic model can be discussed as $V^{-1}UV=W$, with $dW \notin G_{H}$, for instance $dW - dI \in G_{H}$.
ex: Assume existence of pseudo base, analytic over $\Omega$, further F a linear transformation in $\gamma$, $F(\gamma)(\zeta)=f(\zeta)$ and choose $d U \in G_{H}$, with $\gamma(u)$ analytic over $\Omega$.

\section{The polar set}

Assume $f \in (I)_{\mathcal{H}}$ an ideal with topology $\mathcal{H}$ and $d U \in G_{\mathcal{H}}$.  Assume $\Phi$ reciprocal linear transformations over G, $\Phi(G)=F$ and $\Phi(G')=F'$.
Relative a decomposable topology, $\mathcal{H}(G \times G') \simeq \mathcal{H}(G) \times \mathcal{H}(G')$,
$(G,G')$ can be determined by $\Phi^{-1}(F \times F')$. ex $\Phi(dU)=\rho dU$, with $\rho,1/\rho \in L^{1}$.
Given $F=\Phi(G_{2})$ and $F \rightarrow F^{\bot}$ reciprocal, where
$dI \in G_{2}$ with $F \cap F^{\bot}=\{ 0 \}$ and $F^{\bot}=\Phi(G_{r}')$. With these conditions, we have $G_{2} \bot G_{r}'$, in the sense that $< U_{2}f,U_{r}' g>=< F f, F^{\bot} g>=0$.

Consider $dI \in H(C_{K})$, where $C_{K}$ upper analytic envelop (\cite{Bruhat57}), $C_{K}(F,E)=\cap_{E \backslash F \cap K \subset A} A$,
where A analytic and K a compact set.  ex :$d I \notin H(F^{\bot})$, but $dI \in H(C_{K})$.
Thus, an approximation property for $C_{K}$, does not imply an approximation property for $F^{\bot}$. f polynomial on $F \cup C_{K}$, does not imply f polynomial over E.  ex : $\int_{\Omega} P dx=0$ implies $m \Omega=0$. When $\Omega=F \cup F^{\bot}$ and $\Omega'=F \cup C_{K}$, we have $\int_{\Omega'} P dx = \int_{\Omega} \tilde{P}$, where $\tilde{P}$ is not necessarily a polynomial.

ex : $rad (I)(F) = \{ f^{N} \mbox{ reduced } \}$, where $F$ an algebraic set, in this case $C_{K}(F,E)$ is sufficient to define rad (I).
The codimension can be determined for $C_{K}$, that is over regular approximations.

Assume f(x,y,z)=c has $x^{'2}(s) + y^{'2}(s) + z^{'2}(s)=1$.
ex : $\frac{d}{dt} \Sigma x_{j}^{2}=1$ iff $2 \Sigma x_{j} x_{j}'=1$.
When $\frac{d}{dt} < x,x^{*}>$ is (separately) absolute continuous, we have $\frac{d^{2}}{d t^{2}} < x,x^{*}>=0$ implies $\frac{d}{dt}<x,x^{*}>=1$.
Given $\frac{d}{dt}<x,x^{*}>$ single valued in AB (analytic and bounded) in $(x,\frac{d x^{*}}{dt})$ and maps $0 \rightarrow 0$, we have  $\frac{d}{dt}<x,x^{*}>$ is linear in
$(x,\frac{d x^{*}}{dt})$, that is $\frac{d}{dt}<x,x^{*}> \sim <x,\frac{d x^{*}}{dt}>$. Further, given canonical dual
$<\frac{dx}{dt},x^{*}> \sim < x,\frac{d x^{*}}{dt}>$. We have $\frac{dx}{dt} \bot \frac{d x^{*}}{dt}$ implies
$\frac{d^{2}}{d t^{2}} < x,x^{*}>=0$ and given $x^{*}(dx/dt)$ absolute continuous (and $x(d x^{*}/dt)$ absolute continuous), we have $\frac{d}{dt}<x,x^{*}>=1$.

Starting from $V(f)=\{ p(x,y)=0 \}$, assume $\eta$ a formal analytic solution to the algebraic equation p=0. Then, $\eta(t,y)=a_{0}(y) + a_{1}(y) t + a_{2}(y) t^{2} + \ldots$. When $a_{1}(y) \not\equiv 0$, then V(f) has property (H) (\cite{Oka60}).  ex : $\eta$ analytic in t=$x^{1/\mu}$.

Ex : Assume for $f,\rho$ BV, $J(f)=\{ \rho \quad \rho f \mbox{ BV } \}$, $\frac{d (\rho f)}{\rho f}= d \log \rho + d \log f$. Let $M(f)=\sup_{\Omega} \mid f \mid$, that is $\mid \int d \rho f  \mid \leq M(f) T(\rho)$. Further, $T(\rho f) \leq M(f)T(\rho) + M(\rho) T(f)$. In particular, $\frac{T(\rho f)}{M(\rho) M(f)} \leq \frac{T(\rho)}{M(\rho)} + \frac{T(f)}{M(f)}$.
Note that $\log f \in L^{1}$ implies $f \in \mathcal{D}_{L^{1}}'$.

We will see that characteristic for a HE (hypoelliptic) symbol f, is monogenity in the representation space. Regularity is preserved under very regular change of local coordinates (\cite{Dahn24}). Note that when G is continuous over H, $dI \in G$ has resolution in 1-parameter subgroups. Monogenity (a contraction property) implies separate continuity over 1-parameter subgroups.

Assume $\mathcal{F}$ analytic and single valued functions (\cite{Ahlfors_Beurling}). The class is called compact, if given open sets $\Omega_{n} \uparrow \Omega$, and given a uniformly convergent sequence relative $\mathcal{F}(\Omega_{n})$, we have existence of $g \in \mathcal{F}(\Omega)$, so that $\max_{\Omega} \mid g(x) \mid=\mid f'(z_{0}) \mid$.

Assume $M_{\mathcal{F}}(z_{0},\Omega)=\sup_{f \in \mathcal{F}} \mid f'(z_{0}) \mid$. $N_{\mathcal{F}}$ defines a class of zero sets, where $M_{\mathcal{F}} \equiv 0$ on $\Omega$. D denotes the class of analytic functions with finite Dirichlet integral. AB denotes the class of analytic functions bounded on $\Omega$.

\newtheorem{lemma10}{Lemma}[section]
\begin{lemma10}
For a homogeneous domain E without irreducibles, we have $E \notin N_{D}$, that is $M_{D} \neq 0$ on E.
\end{lemma10}

Consider a homogeneous space (x,y(x)) and $M_{D}(p,\Omega)=\sup_{f \in D(\Omega)} \mid f'(p) \mid$.
Assume for instance E constant surfaces to $f \in D$.
Hurwitz gives as $M_{D}$ uniform limes on compact sets, that $M_{D} \equiv 0$ or $M_{D} \neq 0$.
Assume f($G_{2}$)=const of maximal order and $\frac{df}{dn}=\Sigma \frac{\delta f}{\delta u_{k}} \frac{\delta u_{k}}{\delta n}=0$ with $u_{k}$ regularly dependent of n, that is $\frac{\delta u_{k}}{\delta n} \neq 0$. Given $M_{D}(dN)=0$ we must have $ dU_{j}/dN = const$.

Annihilator representations, are linear representations. Continuation using a linear representation in phase is corresponding to algebraic continuation.  $Uf = \int f(x) d \alpha(x)=f(\xi) \int d \alpha(x)$ $(+ \int df $ flux), for some $\xi$ is corresponding to an interpolation property.
Using algebraic geometry $I_{AB} \subset I_{D}$ implies $Z(I_{D}) \subset Z(I_{AB})$. Generally IZI=rad (I).
For linear sets we have $M_{AB}=M_{D}=0$ or both $\neq 0$.

Assume $\chi_{E}$ characteristic for E, with $E=\{ \mu < f(\xi) < \lambda \} \subset \subset \mathbf{R}^{n}$. Further, consider a set F with $\int \int \chi_{E} d \mu d t/t$, for $d \mu$ Borel, with $\chi_{F}=\chi_{E}/t \rightarrow 0$, as $t \rightarrow \infty$.
Necessary for $I(F) \subset I(E)$ is a strict condition, for instance algebraic polar. Further, when E is in a normal rectifiable set, I(E) can be represented with property (H) relative F. ex : Assume $B_{c}=\{ \xi \quad \mid \xi \mid^{c} < R \}$ and that $\{\xi \quad \mu < \mid f(\xi) \mid < \lambda  \} \subset B_{c}$, for some positive c,R. This condition is satisfied by HE symbols. When E is compact, it can be proved (\cite{Nils80}) that bd E is (n-1) dimensional.
Note that $\chi_{E}/t \rightarrow 0$ does not imply a rectifiable boundary.
However, assume $\chi_{E}=t \chi_{F}$, when $g_{E}=\mathcal{F}^{-1} \chi_{E}$, we consider $\chi_{E}$ as $\mathcal{F}(\frac{d}{dn} g_{F})$. Given monogenity, we can prove $< \chi_{F}(f),dN>$ is regular, where $dN>0$ and f regular.

\section{Linearity}
Linearity in $I_{AB}$ does not imply linearity in $I_{D}$. $M_{AB}=0$ implies $M_{D}=0$,
$I_{AB}=\{ f \in (I) \quad M_{AB}=0 \mbox{ on } \Omega \}$. When f AB is single valued outside $\Omega$, then f is linear. Further, $I_{AB}(\Omega) \subset I_{D}(\Omega)$ .

\newtheorem{lemma12}[lemma10]{Lemma}
\begin{lemma12}
Conjugation relative a bilinear form, according to $dU+dV=dI$ locally, generates a triangulable representation space.
\end{lemma12}

Starting from $0 \leq <x_{j},x^{*}_{j}>=\Sigma \lambda_{j}^{2} \leq 1$, with canonical dual topology.
Assume dU,dV conjugated according to a bilinear form,
further volume preserving relative scaling, in particular the conjugation maps lines on lines. Then the representation space can be seen as an abstract complex (\cite{AhlforsSario60}). Assume $dU+dV=dI$ locally. When $G_{2}$ continuous, the representation space is triangulated (\cite{AhlforsSario60}).
When, h : $G_{2} \rightarrow (x,y,z)$ (\cite{BrelotChoquet51}) is a homeomorphism, giving a very regular domain, the triangulation of representation space induces a triangulation (stratification) in 3-space (\cite{AhlforsSario60}). We can (\cite{RieszM37} extend to G, that is $\Sigma \lambda_{j} M_{j}$, where $M_{j}$ are defined by sub groups to G.

Assume sng f=$\{ f=c \}$ is given an absolute continuous representation according to $\{ df =0 \}$. Given $\Phi$ euclidean, with $\Phi(f) =\int f dU= c \int dU$, we can determine $\Omega$, with $\int_{\Omega} dU=1$, a ``normalized'' surface.

E is called K-irreducible, when $E \subset F \cup F'$ implies $E \subset F$ or $E \subset F'$,
(\cite{Bruhat57}).
 ex : $\{ \mu < f < \lambda \}$, that is when $F=\{ f=\mu \}$,
F is a K-irreducible composant in an analytic set E, implies $F \not\subset C_{K}$.
Presence of trace implies the set is not K-irreducible.  ex: C cylinder web in representation space, identifies inner and outer domains (locally), but the disk closure of $C \backslash \{ u=v \}$ does not exclude spirals.
Given f is regular,
$\{ f=c \}$ has measure zero (Sard).  f BV, does not imply that f is regular, why
f=c is possibly non-trivial.

The trace $\Gamma$ is seen as irreducible. Given $\Gamma$ a 1-dimensional generatrix, relative a continuous group G, a necessary condition for $h(\Gamma)$ to be a 2-dimensional surface, is free movement in G, for instance a surface of revolution. ex : $\Phi(dU,dV)=$
$(dU,dV^{-1} U)$, that is $d V^{-1}U \rightarrow dI$ generates $\Gamma$.

Multivalentness with holomorphic leaves have spiral approximations.
$C_{K}$ according to analytic disks, does not exclude a trace in the polar.  Monodromy implies $(\gamma,f(\gamma)) \rightarrow \gamma$ single valued. HE implies monodromy for spectral resolution.

Assume $T(f)=\int \mid f'(z) \mid \mid dz \mid $ and $D(f)=\int \mid f'(z) \mid^{2} dxdy$.
bd U=$\{ D=\mu \}$ and bd $U'=\{ T=\mu \}$, but bd $U' \not\subset$ bd U or may have no points in common.

\newtheorem{lemma11}[lemma10]{Lemma}
\begin{lemma11}
When f is algebraic in (u,v) conjugated, we have existence of path in $bd U'$ $\cap$ bd U.
\end{lemma11}

We will later show that monogenity is sufficient for a non-trivial intersection. Linear phase implies algebraic action. More precisely, when U has algebraic action in the sense that UI=IU relative the (Laplace)
transform $\widehat{I}(\phi)=e^{\phi}$, Then, U can be determined using $U e^{\phi} \simeq e^{U \phi}$.
When f is regular in a neighborhood of UI=IU, that is f(u,v) is regular, when $v \rightarrow 0$, then f can be seen as monogenic.
ex : $e^{P}=const$ implies $P'=0$, given P absolute continuous, thus P=const.

\newtheorem{ex1}[lemma10]{Example}
\begin{ex1}
Rectifiable parabolic sets are not symmetric.
\end{ex1}

Rectifiable curves are not necessarily continuous. Assume $\phi$ are BV, then $\gamma=(\phi,0,0)$ are rectifiable and $T(\gamma)=s$ (arc length) (\cite{Riesz56}). When $\gamma$ defines a parabolic surface, with $\phi$ SH (subharmonic), we have
$\phi<0$ implies $\phi=const$. Further, we have presence of a maximum-principle. Close to maximum points, we do not have symmetry for the domain of curves.

\section{An interpolation property}

We define $T_{A}x=\lim_{s \rightarrow 0} \frac{1}{s} (T(a(s)) -1)f$, where $A=G_{1} \subset G$ (\cite{Garding47}).
For a 1-parameter group A, tangents are given by $dU(f)=0$, $d U \in A$. We have an
inversion formula, $f(h)=\int h(b) T(b) f d b$, where db volume on G, left invariant, $h \in C_{0}^{r}$
and $f \in B$, Banach-space.

Assume db volume on G, corresponding to conjugation in G. Necessary for a well defined inverse, is that tangents through a homeomorphism (\cite{BrelotChoquet51}) correspond to a 2-dimensional surface. ex : db very regular subgroup (\cite{Dahn24}).

\newtheorem{lemma3}[lemma10]{Lemma}
\begin{lemma3}
Given (dU,dV) BV measures, with multipliers in $L^{1}$ and with an interpolation property in $G_{H}$, then a 2-surface is generated in 3-space.
\end{lemma3}

Perfect sets have an interpolation property. Assume existence of dW $\bot dU \times dV$ and $dW \bot dV \times dU$. Further, $dU/dW=\sigma$, $dW/dV=\mu$, with $\sigma,1/\mu, \sigma \mu \in L^{1}$, that is $\rho = \sigma \mu \rightarrow 0$. When $\rho \in H$, without essential sng in $\infty$, we have  an interpolation property in $G_{H}$. Existence of dW non-trivial, implies that (dU,dV) generates a (0-) 1- or 2-surface in 3-space. When both $\sigma,1/\mu \rightarrow 0$, we have a 2-surface in 3-space.

Assume $\gamma(t)$ a trajectory and $t \in \big[ 0,1 \big]$, $\int_{\gamma} f dU(x)=f(\xi) \int_{\gamma} d U(x)=f(\xi)\int_{\tilde{\gamma}} dI(x)$, for some $\xi$, where we can assume $\int_{\tilde{\gamma}} dI=1$.
Thus, $\tilde{\gamma}$, is defined by ${}^{t} U \gamma$ and given $d {}^{t} U$ preserves regularity, the trajectory is connected. Thus, $Uf(\gamma)=f(\xi)$, that is $\xi \in$ multivalentness,
(Abel's problem). Monodromy implies $\xi \in \tilde{\gamma}$. Thus, closedness is a monodromy property.
It is an interpolation property, by a continuous Stieltjes representation.  A very regular subgroup does not have a trace. $\gamma_{1} \sim \gamma_{2}$ refers to existence of a regular deformation, this implies interpolations property. that is a `` relative approximation property'' or
a contraction property.

 Existence of dN normal to (dU,dV), that is simultaneously normal to dU,dV, is implied by monogenity. Given dV/dN has a regular inverse (two mirror model), dN can be seen as an intermediate measure, that is absence of a maximum-principle and presence of an interpolation property in G.  For a rectifiable surface and a curve on the surface, with an interpolation property (absence of trace), the curve is not a geodesic in limes of strata (\cite{ Backlund16}).

We can discuss singularities as relating to type of movement, for instance given bounded singularities relative dU, we have that $dU^{N}$ give regular approximations. ex : When ker $dU^{N}-dI \downarrow \{ 0 \}$, singularities are considered as bounded relative dU.

Dirichlet domains are not perfect. $D_{1}(f)=\int \mid f'(x) \mid^{2} dx$ is corresponding to
$< \delta(x-y), \mid f'(x,y) \mid>$ (\cite{Dahn25}), that is with support on the diagonal. For an extremal definition (point support), we do not assume that $\{ f=c \}$ separates points in the domain.
When $\{ f=c \}$ is rectifiable and n-1 dimensional, it separates points in the n-domain.

For a very regular subgroup and dU+dV=dI, we assume that points are not irregular for dU,dV simultaneously.
That is change of type is necessary according to the two mirror model (\cite{Dahn14}), cf density for $G_{2}$. Monogenity implies contraction to $G_{1}$, under preservation of regularity.
Further, monogenity implies uniform convergence in inner points.

\section{A Brelot mapping}

Assume E a connected topological vector space,  with dimension $\tau \geq 2$. Further, existence of a homeomorphism that maps nbhd $P_{1}$ in $G_{\tau}$ on nbhd $P_{1}'$ in $\mathbf{R}^{\tau}$ (possibly compactified), with compatibility condition $h_{P_{1}} h_{P_{2}}^{-1}(C)$ a local isometry (or when $\tau=2$ locally conformal), where $C = nbhd P_{1} \cap nbhd P_{2}$. E with this structure, is countably compact and separated. A regular boundary, is defined as the image under h, in a $\tau-1$ dimensional surface. A regular surface is a relatively compact open set, one sided relative the boundary, image of a homogeneous space, for instance $z(x,y) \in C^{2}$ defines (x,y,z). Inner points in nbhd P are mapped on inner points in nbhd $P'$
locally 1-1. A regular domain is very regular if also the boundary is regular (\cite{BrelotChoquet51}). A geodesic metric is defined using the distance  $l_{g}=\int_{\gamma} ds$, where $d s^{2}$ gives euclidean distance.

ex : Assume T defines rectifiable trajectories and $T(dV)(p)=T(dV^{\diamondsuit})(\overline{p})$. $w=(p,q)$ and
d is the distance between p,q, $d(p-q) + d(p+q) \leq T(w) + T(w^{\diamondsuit})$. Note that MV(p,q)=$(p+q)/2$, that is MV(p,q) gives lower bound, when T(f) is HE and p=q.

Assume $(I)(dU,dV)$ with a decomposable topology, further $UVe^{\phi}=e^{(U+V) \phi}$, that is $\phi \in (I)(dU,dV)$ implies $e^{\phi} \in (I)(dU,dV)$. This defines a convolution algebra, when the topology for (I) is dual to $\mathcal{E}^{'(0)}$.
When $U \widehat{I}(\phi)=\widehat{I}(U \phi)$, that is when U does not change type, the entire space is not reached and we have a non-trivial polar.  When $U \widehat{I}=\widehat{I}V$ and (dU,dV)=(dV,dU), the polar can be given in Exp.

\section{Rectifiable sets}

$\Gamma= \{ (x,y,z) \}$  is rectifiable, if $x^{'2}(s) + y^{'2}(s) + z^{'2}(s)=1$. $T(f)=\int \mid df \mid$ defines a rectifiable  curve $(f,0,0)$. ex : When A is a rectifiable 2-dimensional surface in 3-space, A can be generated by its tangent planes.
A rectifiable implies that A separates points in $\mathbf{R}^{3} \backslash A$.

Assume existence of h  homeomorphism (\cite{BrelotChoquet51}) $h(G)=V^{n-1}$, such that a regular boundary is mapped on a (n-1)-dimensional set in n-space. Assume $0 \in V^{n-1}$ and $V^{n-1}$ symmetric. Given $V^{n-1}$ is simply connected, by Riemann mapping theorem : $V^{n-1} \simeq B(x,t)$ a disk, using a holomorphic  mapping.

Let $\big[ d U,d V \big] f=\frac{d^{2} f}{d t^{2}}$ and $d Vf=\frac{df}{dt}$. We then have (\cite{Lie91}, chapter 18)
$\big[ d U,d V \big]= d V$ or =0. In the first case for instance $\frac{df}{dt}$ absolute continuous. $\varphi_{1} dU_{1} + \varphi_{2} dU_{2}=\omega dI$
means $\omega=0$ or $\omega \neq 0$, for instance $\omega \equiv 1$. Given $\varphi_{j} \in L^{1}$, we have a local deformation, where the movements have different trajectories.

We assume rectifiable surfaces (curves) $\{ f_{i}(x,y)=c \}$, are such that $df_{i} \neq 0$ on x=y, that is evaluation is according to $< H(x) \otimes H(y), df>$.
Assume dx/dt=X, dy/dt=Y and (\cite{Dahn13}) consider the projection method. When f HE, all iterated arithmetic means are of real type and we have a stratification of real type.
Consider an absolute continuous deformation according to $d I \rightarrow \rho dI$ with $\rho \in L^{1}_{ac}$. Given a closed curve, we have when
$\int_{\Gamma} \rho dI=\int_{\Gamma} \rho'(x) dx=0$, defines an exact form, that is $\rho'=0$ over $\Gamma$, implies $\rho=0$. Rectifiable surfaces (curves) have absolute continuous deformations (\cite{Backlund16}).
In particular, f HE implies that singularities can be represented on rectifiable surfaces.

A rectifiable  curve has
$\mid I_{j} \mid < \mid \gamma \mid$, where $I_{j}$ are chord segments, that is $\gamma$ arcwise convex (concave).
Assume $\tilde{\Gamma}$ the smallest closed convex polygon, that includes vertices on $\gamma$,
further that $\gamma$ has positive curvature and
$\frac{1}{\mid \gamma \mid} \int_{\gamma} = \frac{1}{\mid \Sigma I_{j} \mid} \int_{\cup I_{j}}$, then we have
$(\cup I_{j}) \subset (\gamma)$.

$Uf=\int f dU + \int_{\Gamma} df$, $Uf \rightarrow f$ regularly, assumes a rectifiable boundary.
f absolute continuous over $\Gamma$ implies $\int_{\Gamma} df=f$. When $\Gamma$ multivalent,  the  curvature is multivalent. The total curvature, is the sum of curvatures over segments. ex : $\mid \rho \mid=c_{1} \rightarrow \mid \rho \mid=c_{2}$, a movement along the normal. A sufficient condition, for the normal to have a unique orientation, is monogenity for the normal operator. Sufficient for this, is a regular polar.

Assume S a characteristic surface and consider sections according to $\mathcal{F}(x_{0})=\{ y \quad (x_{0},y) \in S \}$. Further, assume the sections bounded. When $\omega$ dicylindrical is of class (H) (\cite{Oka60}), we have that  log diam $\mathcal{F}(x)$  is locally SH in x.

\newtheorem{lemma15}[lemma10]{Lemma}
\begin{lemma15}
Assume f regular analytic and  $M_{\lambda}=\{ \mu \leq f(x,y) \leq \lambda \}$
compact, then bd $M_{\lambda}$ can be given by rectifiable and normal curves.
\end{lemma15}

The result is a direct consequence of (\cite{Nils80}).
As f HE, we have for
f=const, that $If=\int f dI + \int_{\Gamma} df=const$. Assume $\Gamma$ defined by g regular, with $\{ g,f \}=0$,  for instance Xdy-Ydx=0, that is df=0 on $\Gamma$. Thus, relative g,  $\Gamma$  is rectifiable.
It follows given f HE, that the rectifiable first surfaces are normal. In this setting, $\Gamma$ rectifiable, is dependent of type of dU, $dU \in G_{H}$ continuous.

\section{Normal elements}
A normal rectifiable family $\{ \gamma \}$,  is such that $T(\gamma)< C$ implies $\mid \gamma \mid < C'$, or
every infinite sequence $\gamma_{n}$ has a limit, uniformly continuous in inner points (\cite{Montel24}). ex : $f=\rho f'$, with $\frac{d}{dt} \log f=1/\rho$. Let $\rho \rightarrow \infty$, then when f rectifiable, it is not necessarily normal.

D(f) has extremal definition in the following sense, $\mid f'(z) \mid \mid dx \mid \simeq  d \mu_{x}$, means that $< \delta(x-y), d \mu_{x} \times d \mu_{y} >=\int \mid f'(z) \mid^{2} \mid dx \mid$. N is referred to as perfect, if $<N,d \mu_{x} \times d \mu_{y}>=0$, implies
$< N,d \mu_{x}>=0$, $<N,d \mu_{y}>=0$. Assume $\int_{\Gamma(x)} \rightarrow \int_{\Gamma(y)}$ projective. When f HE and absolute continuous, we have $df=0$ on $\Gamma$ implies $\Gamma$ rectifiable, when the condition is symmetric, we have that $\Gamma(x,y)$ is rectifiable.

Assume $D(f)=\int \mid X \mid^{2} + \mid Y \mid^{2} d x dy < \infty$. Given $\dot{x}^{2} + \dot{y}^{2}=1$
implies $(x,y) \in B$, this defines a normal rectifiable family. As $\dot{x}^{2} + \dot{y}^{2} \rightarrow 0$, we do not have $(x,y) \in B$, that is
a normal surface does not define a normal rectifiable family.
Note that $(x,y,f(x,y))$ is rectifiable, where f(x,y)=const,
only when also $(x,y)$ rectifiable.

Assume S first surfaces to an holomorphic function F.
For a normal family of curves, we have that the total length of curves is bounded. F HE implies that F is unbounded, that is the length has no relative maxima in the finite space. Thus for HE symbols, we have that dist($\{ f=c \}$,0) is finite, for finite c.

ex: Determine Q reduced, such that f/Q=g is BV, that is rectifiable. Thus, invariants over f  are invariants over g, that is $\Sigma(g)=\Sigma(f)$. $\mid f \mid < C$ implies $\mid g \mid < C'$, for $\mid x \mid >R$. $1/Q \rightarrow 0$, corresponding to an inclusion condition, that is $g \rightarrow f$ projective. Assume $f(u_{1},\ldots,u_{p})=0$. Presence of a chain
$u_{1}=0 \rightarrow u_{2}=0 \rightarrow \ldots \rightarrow u_{p}=0$, is dependent of the approximation property,
that is a resolution of zeros, corresponding to resolution of dI. A resolution that preserves monogenity assumes an initial condition, corresponding to $f(u_{1},\ldots,u_{p})=0$, $\frac{\delta f}{\delta u_{1}}=0$ implies $\frac{\delta f}{\delta u_{j}} \neq 0$, for some j. Further, that the chain is absolute continuous, more precisely when $dU_{2}=\rho d U_{1}$, $\rho \neq$ constant.

ex : $F(\gamma) \rightarrow I$ defines a neighborhood
of $F(\gamma)=I$. When the approximation is not regular, there is room for spiral approximations. Presence of trace in representation plane, is indicated by geodesic  curvature =0. When $F(\gamma) \rightarrow I$ is BV, geodesic curvature is invariant and a necessary condition for regular inversion, is absence of trace.

\newtheorem{lemma7}[lemma10]{Lemma}
\begin{lemma7}
D(f) defines a semi-norm.
Regularity for dD(f) does not imply regularity for dT(f) close to first surfaces. Assume $dD \rightarrow dT$ projective and that $\{ D=c \} \rightarrow \{ T=c \}$ absolute continuous (non-constant), then monogenity is preserved close to first surfaces to T.
\end{lemma7}

According to en earlier lemma, when f algebraic in (u,v), monogenity is preserved. Assume dD monogenous and that we have existence of $dN \bot dD$, which implies $dN \bot dT$.
In particular $<N(f) , dD>=0$ implies $<N(f), dT >=0$ close to $\{ D=c \}$.
A monogenous function in composition with an absolute continuous mapping (non-constant), is monogenous.

Fix p+1 points, $\frac{d^{m} \Phi}{d t^{m}}f=c_{m}$, for constants $c_{m}$.
$\Phi^{(p)}f \neq 0$, with p+1 initial conditions, define a normal family (\cite{Montel_37}). Evaluation does not identify invariants, we will argue that a sufficient condition, to identify movements in a stratification,  is monodromy. Monodromy implies a regular contraction. Abscense of monodromy for a (spectral) stratification, implies a not HE symbol. HE refers to a reduced symbol
over a domain of order 0, for instance a disk domain. Order for a group is dependent of topology, in $L^{1}$ we can consider $G_{2}$, why a second order condition  is sufficient, to conclude HE.
p+1 linearly independent movements define a continuous group, that is an approximation property. This condition does not imply HE.

A representation in the plane corresponding to $G_{8}$, defines a quasi-normal family (\cite{Montel24}). Contraction to $G_{2}$ is corresponding to uniform convergence in higher order derivatives. $dU=\rho dI$, with $\rho \in L^{1}$ is approximated by H.
Inner points correspond to a strict condition, which gives a contraction property.

$e^{-px} T_{x} \in \mathcal{S}_{x}'$,  is $C^{\infty}$ in $\xi$, as $\xi \in cl(\xi_{j})_{1}^{p+1}$ (convex cone), with $\xi_{j} \in \Gamma$, for a $\Gamma \subset \mathbf{R}^{n}$ (\cite{Schwartz52}), that is $\xi$  is the result of a local change of local variables $\xi_{j}=\Psi_{j}(x,y)$ and $p=\xi + i \eta$. As $\xi \in$ inner points in $\Gamma$, we have $e^{-px}T \in H(p)$. p such that $e^{-px} T_{x} \in \mathcal{S}_{x}'$ defines a cylinder $\Gamma + i \Theta$. Assume $G_{2}$ according to the above, with conjugation according to $\Psi_{j}U=V \Psi_{j}$. Under the condition for continuation using dU, $\mid \tilde{f}(\xi + i \eta) \mid < P(\eta)$, for P polynomial, we can determine continuation by dV, according to the above.

\newtheorem{lemma17}[lemma10]{Proposition}
\begin{lemma17}
f HE implies f monogenous over $G_{2}$, given $G_{2}$ very regular and abelian.
\end{lemma17}

Given $G_{2}$ abelian and very regular (\cite{Dahn24}), HE is preserved. Further, there is dN with $dV \prec dN \prec dU$, that is we can assume $dV/dU \rightarrow 0$ regularily in $L^{1}$. Thus f is HE over (dU,dN) and (dN,dV).
We have separate regularity, when we let $dN \rightarrow dI$ regularily. Further,
$\frac{\delta f}{\delta n}=\frac{\delta f}{\delta u} \frac{d u}{d n} + \frac{\delta f}{\delta v} \frac{dv }{d n}$ is finite. As f HE and (u,v) very regular, we have $\frac{\delta f}{\delta u} \neq 0$ and $\frac{\delta f}{\delta v} \neq 0$ regularly and bounded. As $G_{2}$ abelian, (u,v) defines disk domains, thus the derivate is single valued. cf P-convexity according to (\cite{Dahn22}).

\section{Continuable elements}

We can define subordinate elements, according to $f \prec g$, if $f=g \circ v$ with $\parallel v \parallel < 1$, for instance
existence of $d U \in G$, a given continuous group, with $Ug=f$,  that is g continuable by G.
When f is holomorphic, there is a continuous deformation$,\{ f=c \} \rightarrow \{ f=0 \}$.

Subordinate elements define continuable elements. For our applications, a continuous inclusion $I(dU) \subset I(dV)$ is necessary for continuable elements.
Assume $\{ \Gamma_{i} \}$ defines a very regular boundary, according to $dU_{i}=dI$ and $\{ f=c \}$. Further, assume $\Gamma_{i}$ rectifiable, when $i \neq 1$, thus $\Gamma_{1} \rightarrow \Gamma_{i}$ absolute continuous. Given all $\Gamma_{i}$ are rectifiable, we have reversible continouability.
Existence of $G_{2}$ very regular, gives continuable elements.

ex: continuability according to divisibility for modules $(b) \subset (a)$.
$d V^{n} =\rho_{n,m} dU^{m}$, with $\rho_{m,n} \in L^{1}$.

ex : $f'(\zeta + e^{-s})$ according to Dirichlet series, with $s=t + i \sigma$, single valued as $t > \rho$, where $\rho$ is radius for convergence.

ex: f(z) absolute continuous or f(1/z) absolute continuous. f(1/z)=1/g(z) absolute continuous means that $dg/g^{2}=0$ implies 1/g const implies g=const, given $g^{2}$ bounded.
When we consider $\{ f(u)=const \} \rightarrow \{ f(v)=const \}$ absolute continuous, where f is absolute continuous in v, note that f analytic in u
does not imply that f is analytic in x.  ex : contraction, $\frac{(\delta f/\delta u)(\delta u/\delta x_{j})}{f(x)} \rightarrow 0$, given $\frac{f'(u) du}{f(x) d x_{j}} < \infty$. This is the case, when f is HE also in x.

ex : $T_{\Omega}(1/f)=0$ implies $m \Omega=0$, given $\mid df \mid /\mid f^{2} \mid$ algebraic, for instance on points where $\mid df \mid$ finite and $1/f^{2}=0$.
ex : $T(f)=\infty$ implies $D_{1}(f)=\infty$, but without algebraicity, we do not have $m \Omega=0$. For instance $T(g)T(f)=1$, where $T'(g)$ algebraic, over $\Omega$.

 f HE implies f unbounded. f absolute continuous and rectifiable implies $f \sim \int df< \infty$, that is $T(f)<\infty$ implies that f HE, has a maximal domain for absolute continuity. Further, that the sub level surfaces to f not HE, have clustersets.

Given a strict condition in the plane, with algebraic polar, there is an analytic representation in $L^{1}$.
$\Omega \ni x \rightarrow 1/x \in \Omega$, implies a meromorphic representation.
A strict condition implies continuability to $\infty$. ex: assume invariants S $\subset$ a domain of holomorphy. In particular, F=0 implies dM=0 implies M=const, defines a rectifiable nbhd of singularities. Sufficient for this, is that $F(u,v)$ is algebraic.
We can discuss an approximation property according to $dI \in G_{H}$ iff $S=\{ F=0 \}_{H}$.
When dM is a total differential and V is polar, we have dM=$dM \mid_{V^{c}}$, that is a total differential  is dependent of representation of polar set.

ex : assume triangulable boundary, according to dU+dV=dI, with $dU,dV \in G_{2}$, for instance an extremal ray. When C a cone = $\cap H$, where H half spaces containing C. Assume H is defined by $L_{j}(x)=L(x_{j}) \geq c_{j}$ $\forall j$, where L gives a separating linear functional. Given a strict condition, for instance starting from a resolution $f=\Sigma f_{j}$, where $\mid f_{j} \mid >0$, as $x \neq y$, then $\Sigma \mid f_{j} \mid=0$ can be used in a Hausdorff condition for uniformity.

ex: assume $\phi : \Gamma_{1} \rightarrow \Gamma_{2}$, with $\phi$ absolute continuous. Assume $df \rightarrow d \tilde{f}$ preserves zero lines in $L^{1}$. $ d \tilde{f}=df(\phi(z))$, with $d \phi/dz=0$ implies $\phi=const$ and $\phi=const$ removable.

Assume $\mid f(z) \mid \leq C (1 + \mid z \mid)^{a}/\mid p(z) \mid^{b}$, $\exists a,b$ and $V(f)$ is given by p(z)=0 (\cite{Nils80}).
Outside V,  $pf \in \mathcal{O}_{M}$ is decomposable, that is the property  is dependent of V.
Consider the representation of f and $M_{\lambda}=\{ \mu \leq f(x) \leq \lambda \}$, for some $\mu$ and a (spectral) projection operator (\cite{Nils72}), $e(\lambda)=\int_{M_{\lambda}} g(x) dx$, $\lambda \in \mathbf{R}$. ex : let $t=\log \lambda$, $e'(\lambda)=\Sigma t^{j} h_{j}(e^{t})$, with $h_{j}$ analytical. Note that sharp fronts analogously with $\widehat{h}_{j}(t) \rightarrow \widehat{h}_{j}(-t)$, is corresponding to $z \rightarrow 1/z$. When regularity is preserved, we must have absence of essential singularities in $\infty$, for $f \in H$, that is f=p/q, for p,q entire functions.

Assume $\gamma$ are closed paths outside V, with single valued analytic branches in a punctured nbhd of $\infty$,
assume $T_{\gamma}$ a regular analytic continuation along $\gamma$ and between analytic branches of $\gamma$ (\cite{Nils80}). From monodromy $T_{\gamma} \rightarrow \gamma$ is 1-1 with respect to
$\gamma$ connected. $T_{\gamma_{r}} f=f$, for $\gamma_{r}$ according to the above, implies $\gamma_{r} \sim \gamma$, that is a locally exact form. ex : Regular forms in the plane, with singularities on the diagonal, have regular continuation along rectifiable trajectories, corresponding to locally exact forms.

Given singularities on a bounded set, $\tau^{k}f$ is regular, where $\tau$ corresponds to translation (oriented). When singularities refer to deviance from H, then singularities are dependent of topology.
When continuation is given by translation, $\tau^{k}-Id$ gives a neighborhood of invariants.
ex: when f is algebraic, we have a regular neighborhood of lineality $(\tau^{k} -Id)f=0$.

For the class of continuations (\cite{Nils80}) $f \in P_{m}(\mathbf{C}^{n+l})$ and $g(y)=\int_{E_{n}} f(x,y)dx$, over $E_{n}=\{ (x_{1},\ldots,x_{n}) \quad x_{j} \geq 0 \quad \Sigma x_{j} \leq 1 \}$, we have $g \in P_{m+n}(\mathbf{C}^{l})$. ex: $<g(\phi), dV^{m}>(\psi)=<f,dI \otimes d V^{m}>(\psi \otimes \phi)$.
ex: $f \in L^{1}(dU^{k},dV^{m})$ and $dV^{m} \rightarrow dI$ regularly, as continuation of $dU^{k}$.
ex : $f \in P_{m}(G_{l+2})$ implies $g \in P_{m+l}(G_{2})$.

\section{Perfect transformations}

We say that $dN \bot dV$, if $dV=\rho dN$ with $\rho \in L^{1}$ (pseudo-orthogonal), that is dV is absolute continuous as functional with respect to dN.
dN is said to be perfect, if $dN \bot d \mu_{x} \times d \mu_{y}$ implies $dN \bot d \mu_{x}$ and $dN \bot d \mu_{y}$. Given monogenity, $\frac{d f}{d n}$ is single valued, finite and defines a perfect normal.
Assume $dN \rightarrow dI$, then we have $N(\gamma) \rightarrow \gamma(0) + \int_{\Gamma} d \gamma$.
We can assume for instance $\Gamma=\{ dN=dI \}$.
ex : $d V^{\bot} \bot d U_{k}$, $\forall k$ implies $d V^{\bot}=0$ (implies dV=dI), that is $dV$ is not perfect. Note that
$dV^{\bot}=0$ implies dV=dI is dependent of two scalar products, $dV + d V^{\bot}=dI$ and $dU,dV$ are conjugated.

Starting from convex domains, consider sub modules according to $M=G_{r}$.  Form $(e_{k})=\{ u \in M \quad L_{v_{j}}(u)=0 \quad j=1,\ldots,k \}$ (\cite{Dahn22}), where $v_{j} \in G \backslash G_{r}$ and $L_{v_{j}}$ are scalar products.
We can determine a module $\subset (e_{k})'=\{ L_{v}(u) \quad u \in (e_{k}) \}$, minimal for convexity (\cite{ RieszM37}).
Further, starting from an ideal (I), with decomposable topology, for instance $C^{\infty}_{c}$ or $L^{p}$ over open sets, consider $dU$, such that $\int f dU^{N}=U^{N} f \in (I)$. In particular, we can define $UT(f)$ with $T(f) \in C^{\infty}$. Thus, we have existence of $f_{1},\ldots,f_{r}$, so that
$f(u_{1},\ldots,u_{r}) =f_{1}(u_{1}) \ldots f_{r}(u_{r})$. Consider $M^{-1}=\{ dU \in M \quad Uf \in (I) \}$
and K another sub module. Then, assuming $<dU,dN>(f)=0$ has a regular analytic neighborhood for $f \in (I)$, and in the same manner for $dV \in K$, then over (I), we have $(M + K)^{-1}=M^{-1} \cap K^{-1}$ (\cite{RieszM37}). In particular, we have monogenity in the sense that if given $dU \in M$ and $dV \in K$,  with $dN \bot d U \times dV$, we have that $dN \bot d U$,$d N \bot d V$.

Assume dU+$d U^{\bot}$ =dI + dW, when $dU^{\bot \bot}=dU$, we have $dW=dW^{\bot}$. Presence of trace in the polar is interpreted as non removable polar. ex : $dV=\rho dU$, gives absence of projectivity $dV \rightarrow dU$, where $\rho=0$, given dU is BV.
$dV=\rho dU$ with $dU \rightarrow dI$,
that is $\rho=const$ allows invariants of higher order. $dV^{-1}U=dI =\rho^{-1} \vartheta dI$. $\vartheta \in L^{1}$ implies $\vartheta = \vartheta'=0$ isolated points. Thus, $d V U V^{-1}$ has isolated singularities.

$<f(\gamma),dV>=<f,dV>(\gamma)$, when $\gamma$ regular, assumes $<f,dV> \in \mathcal{E}'$. An approximation principle, assumes a sub nuclear id-transformation, corresponding to an inequality of order 0, that is we assume $dV \in \mathcal{E}^{'(0)}$. For a pseudo convex set, the boundary is of order 0. For our application this means that it is sufficient to consider a 2-sub group, that generates disk subsets.

\newtheorem{lemma1}{Lemma}[section]
\begin{lemma1}
Assume $dU \rightarrow dV$ a
conjugation compatible with perfect transformations, in the sense that $dU/dN=\rho_{U}$ and $dV/dN=\rho_{V}$, for $\rho_{U} \in L^{1}$ and $\rho_{V} \in L^{1}$, analogous with continuability in representation space.
The conjugation induces an interpolation property in representation space, as long as further $1/\rho_{V} \in L^{1}$.
\end{lemma1}

In absence of essential singularities in $\infty$, $\rho \in L^{1}$ can be approximated by $\alpha/\beta$, with $\alpha,\beta \in H$.
When $dN \bot dU \times dV$ according to the above, with $\rho_{U}/\rho_{V} \in L^{1}$, we have $dU \prec dN \prec dV$.
When further $dN \bot dV \times dU$, we must have simultaneously $\rho_{V} / \rho_{U} \in L^{1}$. An interpolation property thus assumes further conditions. When $f \in L^{1}(dI)$, we have
$\rho (1/\rho) f \in L^{1}(dI)$. Note that $dU U^{-1} =\rho (1/\rho) dI$, that is seen as translation. Thus, for monogenity, we have to consider approximative solutions (a deformation), where $\rho (1/\rho) -1$ regular. Note that
Id : $\mathcal{D}' \rightarrow \mathcal{D}'$  (pseudo locality), but $Id - \gamma$, with $\gamma \in C^{\infty}$, approximates $dI \in \mathcal{D'}^{F}$ (for instance $\mathcal{D}_{L^{1}}')$.
When $dU^{*}=dV$, we assume $dU^{*} \prec dN \prec dU^{**}$, for the second equation. Reciprocity assumes $dU^{**}/dU \in L^{1}$.

Assume $(a) \cap (b)$ according to $\rho,1/\mu \in L^{1}$, defines monogenous continuation according to
$(M+N)^{\bot}=M^{\bot} \cap N^{\bot}$, that is $v_{j} \in N^{\bot}$ means $L_{v_{j}}(b)$ analytical. If further $1/\rho, \mu \in L^{1}$, we have reciprocity according to
$dU/ dV \rightarrow dV/dU$ in $L^{1}$. Continuation according to $e' \subset e$, when two module elements have a common linear factor, implies an euclidean continuable element. When $\not\exists dW$, with
$dU \prec dW \prec dV$, we have maximality, that is a boundary for continouability (max-principle).

\newtheorem{ex4}[lemma10]{Example}
\begin{ex4}
For every $h \in M \cap N$, we have existence of $n \in N$,$m \in M$, such that 1/m+1/n=1/h.
\end{ex4}

ex : $\{ f=c \}$ defines a normal tube, when the transversal can be reached by points in the first surface.
A normal tube has decomposition in a set with interpolation property (first surface) and isolated points (on the transversal).
Note that a spiral normal assumes a complex representation.

 \newtheorem{def1}[lemma10]{Definition}
 \begin{def1}
 We discuss an approximation property according to a) reduction to a surface of order 0 (cylindrical surface), b) regularity according to $dI \in G_{H}$.
\end{def1}

\section{The contraction property}
As earlier, a family of multivalent functions of order p is called quasi-normal (\cite{Montel24}). It is normal, if we fix p+1 points or points for p derivatives.

When $f'(z) > \mu$ on inner points, this implies $\{ f'(z) \}$ is a normal family (\cite{Montel_37}). For this family, there is analogous with Koebe, a contraction theorem, that is $\frac{1}{\mu} \leq \mid f'(z_{1})/f'(z_{2}) \mid \leq \mu$, for $z_{1},z_{2}$ in inner points.
Assume (dU,dV) defines G according to Radon-Nikodym, that is for instance $dU=\rho dI, dV=\vartheta dI$, with $\rho,\vartheta \in L^{1}$.
Assume existence of dN, with $dN \bot dU \times dV$ implies $dN \bot dU$ and $dN \bot dV$, that is
$dU/dN=\sigma_{1}$ and $dV/dN = \sigma_{2}$, with $\sigma_1,\sigma_{2} \in L^{1}$. Given
$\mu^{-1} < \mid \sigma_{1}/\sigma_{2} \mid < \mu$, we have $dU/dV \in L^{1}$, by contraction, that is we have a homogeneous representation space. Assume further $dN(u,v) \rightarrow dI$ regularly and separately regular, according to Radon-Nikodym,
then $1/\sigma_{1},1/\sigma_{2} \in L^{1}$.

Define for a scalar product $L_{v}(u)$, $M^{-1}=\{ v \quad L_{v}(u) \mbox{ entire regular } \forall u \in M \}$, a closed module.
Assume continuation relative a very regular sub group,
according to $(N+M)^{-1}=M^{-1} \cap N^{-1}$, for instance $(\rho + 1/\rho)^{-1} \in L^{1}$
implies $\rho,1/\rho \in L^{1}$. Given $m_{j}$ a base to M, there is a base to $M^{-1}$, $m_{j}^{*}$ with
$1 \leq \mid m_{j} \mid \mid m_{j}^{*} \mid \leq C$ (\cite{RieszM37}).

Assume $(\delta f(u,v)/ \delta u)/(\delta f(u,v) / \delta v)=dv/du$, where $\{ f' \}$ a normal family, according to Koebe. Thus, when f is symmetric in (u,v), $du \rightarrow dv$ is a normal  mapping,
that is the conjugation is normal. Note that given $d v/d u$ non zero, we may over a planar domain, have dv/du=c with
$\mu < c < \lambda$, that is the condition does not exclude a trace. Given u,v linearly independent, the trace must be bounded (continuum). Given f holomorphic (and $f'$) however, $f' \neq const$.

The condition $\frac{f'(z_{1})}{f'(z_{2})}=\frac{d z_{2}}{d z_{1}} \neq 0$ regular, implies $z_{2}(z_{1}) \in C^{1}$, which implies a homogeneous space. Given dU,dV linearly independent over H, they do not have common invariants, in particular $f'(z_{1})=0$ implies $f'(z_{2}) \neq 0$. ex : $g=\log f$ normal, with $\frac{1}{\mu} < \mid \frac{g'(z_{1})}{g'(z_{2})} \mid < \mu$ and
where $f'(z_{1})=f'(z_{2})$, we have $\frac{1}{\mu} < \mid f(z_{2})/f(z_{1}) \mid < \mu$.

ex : starting from $(u,v,z,\zeta)$, with $z=z(u,v),\zeta=\zeta(u,v)$,
then we have that $dz(u,v)=z_{u}' du + z_{v}' dv$, $d \zeta=\zeta_{u}' du + \zeta_{v}' dv$. When $\frac{1}{\mu} < \mid dv/du \mid < \mu$, the domain where dz=$d \zeta=0$ is a domain where both $z',\zeta'$ are normal.

When $f'(z)>\mu$, we have $<H,f'> \sim < f, \delta>$, that is single valuedness corresponding to
evaluation of f (initial condition).

\newtheorem{lemma4}[lemma10]{Definition}
\begin{lemma4}
We say that the singularities are bounded relative dU, if $d U^{N}f$ is regular, for some iteration index N.
When $\Phi(t) f=\int f d U$, we assume that $\frac{d^{N}\Phi f}{d t^{N}}=\big[ dU,\big[ dU \ldots dU \big] \big] \neq 0$ and regular.
\end{lemma4}

ex : Assume $dU=\rho dI$, with $\rho^{N} \in L^{1}$, for some N. In particular $\mid \rho^{N}  \mid \rightarrow 0$ regularly in $L^{1}$. Note that dU has polynomial coefficients.

ex : $\big[ d U_{1},\ldots,d U_{p} \big] f > 0$, defines a normal model. When we have a contraction to a very regular 2-sub group, we have an approximation property according to the above. Note that global $\not\rightarrow$ normal.

Rectifiability can be related to a linear (disk) condition. Absence of monodromy implies the operator is not HE. Presence of trace (trajectories without dimension for regular approximations) implies absence av monodromy. Assume $G_{2} \subset G$,
where $G_{2}$ includes translation, then an almost periodic (pp) condition, taken as a strict condition is sufficient for density. The dimension is dependent of type av movement!, thus maximal rank is dependent of sub group. Further, bounded singularities are relative sub group.  ex: f=c in 0, $U^{N}f \neq c$.  Monogenity to $G_{1}$ at the boundary is necessary for HE.

Consider $G \rightarrow G_{2}$, with $\int f(u,v) dudv$ finite, $\Phi(U,V)f-\int_{\Gamma} df=\int f dUdV$.
An approximation property for $G_{2}$ determines $f=\lim \int fdUdV=\int f d \Phi(U,V)$. This implies $\Phi$ absolute continuous relative $G_{2}$ (vanishing flux),
as long as $\Gamma$ is rectifiable and $\int_{\Gamma} df$ is bounded.
Assume $Uf=\int f dU + \int_{\Gamma} df$, with $dU=\rho dI$, assume further formally $Uf=\mu f$, that is $d \mu=\rho$. Then $d \rho \neq 0$ for instance $\rho$ regular with $d \rho \geq 0$, defines U as a convex transformation.

Consider $h : G \rightarrow V^{n-1} \rightarrow G_{n-1}$, by for instance monodromy. $dI \in G$ implies
$dI \in G_{n-1}$ given monogenity (regular contraction). Monogenity implies lim $f(u_{1},\ldots,u_{n-1},v)$ = lim $f(u_{1},\ldots,u_{n-1},-v)$,
that is a single valued function (independent of starting point). Obviously $G_{n-1}$ is dependent of choice of $V^{n-1}$.

ex : a generalized cylinder, can be represented by (f,g,n), where $f(x,y),g(x,y)$ are regular curves and rectifiable. In particular when $\delta f/\delta x=0$, $\delta f / \delta y \neq 0$, that is x=y is not in the domain. n is given by lines orthogonal to the curve (f,g,0). The result is a rectifiable surface.

Assume $dU \prec dW \prec dV$, with $dV=\rho dU$, $\rho=\sigma_{1} \sigma_{2} \in L^{1}$ and $dV=\sigma_{1} dW$ and
$dW=\sigma_{2} dU$. The order of a chain, is given by $\rho$ ``irreducible'', that is (dU,dV) a maximal chain and $\rho \neq const$. Thus,  (dU,dV) is not continuable in $L^{1}$.

Assume M a multiplier with 1/M=0 iff $D_{1}(f)=\infty$, for instance $M=\rho + 1/\rho$ multiplier, where $\rho$ satisfies a strict condition. ex : T(f)-$D_{1}$(f) finite, T(f) finite implies $1/M \neq 0$. ex : When T(f) gives arc length, extremal length can be defined by $\sup \frac{T(f)^{2}}{D(f)}$ (\cite{AhlforsSario60}). We have for a curve $\gamma$ with $T(\gamma) < \infty$, that extremal length =0 in $\infty$.
ex : $f=\rho f'$ where $\rho \rightarrow \infty$
($\rho \notin L^{1}$)

Equicontinuous representations over a group $G_{2}$, assumes BV measures. dI is extremal (point support), when $G_{H}$ is continuous. When $dI(x,y)=dI(x) \otimes dI(y)$, it is extremal. In particular, we assume $dI(x,y)=I(dx,dy)$.
Assume $\Phi(f)=\int f dU (+ \int df)$ and existence of $dV \in G$, with $\Phi_{UV}(f)=f(0)$. When evaluation $dI$ extremal (a strict condition), we have density for polynomials in $L^{1}(dI)$. Assume P polynomial, then
PE-I is monogenous when E is parametrix, but not when E is fundamental solution. ex : $PE-I=\gamma$ with
$\widehat{\gamma} \in \mathcal{S}$, when $\gamma=\gamma_{1} \times \gamma_{2}$ and $\gamma_{1} \rightarrow \gamma_{2}$ a linear duality, gives a single valued and finite normal derivative on a cylinder (\cite{Schwartz52}).
When $G_{2}$ is a very regular sub group, $E-I$ is monogenous, if $E-I(G_{2} \rightarrow G_{1})$ preserves regularity.  Assume E a fundamental solution, consider continuation to E parametrix, then for Z, the zero space to E ($\simeq$ polar set), P HE implies Z regular.

Assume Uf regular and bounded in nbhd f=c, that is $df(u,v) \neq 0$ in nbhd singularities to f. A stratification can be constructed using
$\lim_{r \rightarrow \infty} \frac{1}{r} \int_{\Gamma_{r}} df$.  Let $\phi : \Gamma \rightarrow \Gamma^{\diamondsuit}$, that is $d f^{\diamondsuit}=-df \circ \phi$. Assume $\Gamma^{\diamondsuit} \rightarrow \Gamma$ absolute continuous, for instance $\Gamma$ rectifiable, then $\int_{\Gamma} df^{\diamondsuit} \rightarrow \int_{\Gamma^{\diamondsuit}} df$ is projective. ex : dU, $d U^{\diamondsuit}$ closed implies dU harmonic (\cite{AhlforsSario60}) implies MV(f)(u,v)=f(u,v), over real f.

Assume (X,Y) with F Hamiltonian, further according to the projection method (\cite{Dahn13}) $M \rightarrow W$ projective, implies that symmetry for the domain to F (equality for mixed derivatives) is mapped on F harmonic (harmonic conjugation). $W \neq 0$ implies that $(M,W) \neq (0,0)$. A rectifiable stratification starting from the projection method, with for instance $W(\beta f)=M(f)$, assumes $\beta$ non-constant and measurable. For symmetry as dependent of the polar, we assume $dN \bot (M,W)$ implies $dN \bot M, dN \bot W$, that is the normal is perfect as before.
$\widehat{Y}=\eta^{*} \widehat{X}$ gives $\widehat{W}=0$ iff $x^{*} d \eta^{*}/d x^{*}=0$. For preservation of regularity, we assume non-constant multipliers.
$\widehat{Y*X^{-1}}=\eta^{*}$, that is $X \rightarrow Y$ defines a homogeneous space. The condition $\widehat{X*Y*X^{-1}}=\widehat{Y}$ corresponds to conjugation, assuming density for translates.
$M_{1}(X,Y)=H_{x} + G_{y}$, $M_{2}=(\widehat{X})_{x^{*}} + (\widehat{Y})_{y^{*}}$, $M_{0}=(\widehat{H})_{x^{*}} + (\widehat{G})_{y^{*}}$.
$M(X,Y)=M$, $M(H,G)=M_{1}$, $M(\widehat{X},\widehat{Y})=M_{2}$, $M(\widehat{H},\widehat{G})=M_{0}$,
corresponding to resolution of the double laplace transform.

The condition $M_{\lambda}=\{ \mu \leq f(x) \leq \lambda \}$ compact, for some $\mu$, implies exponential type 0 (real type), that is according to the projection method $(X_{0},Y_{0}) \simeq (X,Y)$ (\cite{Dahn13}).
Given conjugation relative Fourier and a HE symbol,
it is to determine a sub group, sufficient to consider $(M_{0},W_{0})$.
The dimension for  $(X,Y) \simeq (X_{0},Y_{0})$ (projection method), gives index for exactness.

\section{Convexity}

f is regular and convex in (u,v), if $f(tu + sv) \leq t f(u) + s f(v)$, with s+t=1. When regularity is preserved, this gives a monogenous convex relation, as $t \rightarrow 0$ and as $s \rightarrow 0$. ex : volume preserving according to $\frac{\delta^{2} f}{\delta u_{t} \delta v_{s}}=\frac{\delta^{2} f}{\delta u \delta v}$.

 We will discuss $f(x,y)$ monogenous in $\omega_{1} \times \omega_{2}$, implies $g(y)=\int f(x,y) dx$ regular in $\omega_{2}$.
ex : $\int f'(x,y) dx=\int f_{y} dI(x)$ and $\int \frac{\delta^{2} f}{\delta x \delta y} dx dy=\int f(x,y) dI(x,y) =If$.

Consider (I)=Exp (I)=$\{ e^{\phi} \quad \phi \in (I) \}$, that is $\phi \in (I)$
iff $\log \phi \in (I)$, that is it is sufficient to consider the phase.
When $< \widehat{I}(g),\widehat{I}(\phi)>=<e^{g},e^{\phi}>$ and when g (positive and real) is convex in $\phi$, we have
$e^{g}=\widehat{F}(\phi)$, with $\widehat{F}$ absolute continuous.
When P is reduced, $\mid e^{\phi} \mid< \lambda$ implies $\mid e^{\phi}/P \mid < \lambda'$, for large x. $e^{\phi}$ convex implies absolute continuous implies outside the polar set, $d \phi=0$ iff $e^{\phi}=\lambda$. Thus, $d \phi \neq 0$ implies $\mid e^{\phi} \mid < \lambda$. $\{ (x,\mu) \quad \phi < \mu \}$ convex implies
$\{ (x,\lambda) \quad e^{\phi} < \lambda \}$ convex.

\newtheorem{prop1}[lemma10]{Proposition}
\begin{prop1}
Assume f monogenous in (x,y,z), then a very regular group induces an absolute continuous mapping to the first surfaces in 3-space.
\end{prop1}

Assume f satisfies the condition of compact sub level surfaces (\cite{Nils80}), then we have that $\{ f= \lambda \}$ is a  compact $C^{\infty}$ manifold, thus rectifiable. Further, we have as before, when f is HE, that f is monogenous over $G_{2}$, very regular and abelian.
$\{ (u,v) \in G_{2} \quad f(u,v)=\lambda \}$ is mapped on a rectifiable surface i (x,y,z).
$f \in C^{\infty}$ has first surfaces of measure zero, by Sards theorem.
Given a very regular representation, every first surface has a regular neighborhood, that is can be considered as rectifiable. More precisely, assume N a normal operator to
the first surface $\Gamma(u,v)$. Monogenity is preserved by a very regular group, why we can assume Nf regular outside $\Gamma(u,v)$.
df=0 on $\Gamma$ as image of $V^{n-1}=\{ f=c \}$, thus $\frac{\delta f}{\delta n} \neq 0$ outside $\Gamma$ , further using monogenity $\frac{\delta f}{\delta u} \neq 0$, as u varies continuously. Thus, for the mapping
$\Gamma(u,v) \rightarrow V^{n-1}$, we see that $f \neq const$ implies $df \neq 0$.

When $\Gamma$ is a very regular boundary, mapped on to a 2-surface $V^{2}$, in particular when we have existence of a 2-sub group, that generates $\Gamma$,
then $d U \rightarrow d U^{\bot}$ must change type. Assume dUV=dVU preserves regularity. BV implies determined tangent. Consider (u,v) conjugated, with determined tangent and such that $G_{2} \rightarrow V^{2}$ preserves determined tangent. Consider h : (du,dv) $\rightarrow$ (dx,dy,dz). When conjugation is given by a harmonic function and h SH, the composition has convex sub level surfaces, which implies a rectifiable boundary.

 f  is SH, given $f=e^{\varphi}$ finite, upward semi-continuous and $\varphi < h(x)$ on bd $\Omega$, with h harmonic, implies $\varphi(x) < h(x)$ on $\Omega$. ex : for line segments $I_{j}$ with end points on $\gamma$,
 assume $\lim_{j} \frac{\mid f(\gamma_{j}) \mid}{\mid \gamma_{j} \mid}=\lim_{j} \frac{\mid f(I_{j}) \mid}{\mid I_{j} \mid}$, where $\gamma_{j}$ is $\gamma$ between the same endpoints,
 then when $\mid f(\gamma_{j}) \mid < \mid f(I_{j}) \mid$ $\forall j$, $\gamma$  is rectifiable.

ex: Consider $e^{\alpha y}u(x)=F_{\alpha}(x,y)$, then we have $\log u$ convex iff $F_{\alpha}(x,x)$ convex. In particular $F_{\alpha}(x,x)$ is absolute continuous.

Assume U defines a Dirichlet domain and K gives a SH continuation to BV, then we have that KU is absolute continuous.
ex : $\int \mid f' \mid (1-\mid f' \mid) dx \leq \int \mid f' \mid \int (1-\mid f' \mid)$, that is given both integrals in the right hand side are finite, we have $D_{1}(f) < \infty$. Assume $\int 1 dx - T(f) < \infty$, then we have that $T(f) < \infty$ implies $\int 1 dx < \infty$ (compact).

\section{Monogenity}
Consider $\int e^{i \tau \varphi}f dx$, with $\varphi \in C^{\infty}$ (\cite{Malgrange73}).
Define $S(\varphi)=\{ H_{1}=\Sigma (\frac{\delta \varphi}{\delta x_{j}})^{2} =0 \}$, for $\varphi \in H(f,S(\varphi))$.
When f is flat on $S(\varphi)$, we can by regular continuation write, $f=\Sigma g_{j} \frac{\delta \varphi}{\delta x_{j}}$, $g_{j} \in C^{\infty}$ flat on $S(\varphi)$ (\cite{Malgrange73}).

When $S(\varphi)=\{ H_{1}=0 \}$, we consider $H_{1} \in H$ as a zero function to $D(\varphi)$, that is on a neighborhood of $S(\varphi)$, $D(\varphi) < \epsilon$.
Given $H_{1}$ polynomial in a neighborhood of $S(\varphi)$ planar, we have that the measure for $S(\varphi)$  is zero.

When $E \subset F$, $dI \notin F'$ but $dI \in E'$, assumes inclusion of the corresponding test space.
The approximation property is sufficient for the implication $dI_{E} \in G \Rightarrow$ $dI_{F} \in G$.

We have seen that a necessary condition for a HE representation, is monogenity over G. Assume $h : G \rightarrow (x,y,z)$.
HE refers to HE over (x,y,z).  HE assumes a contraction property for G, that is G can be reduced to $G_{2}$ translation and
rotation. Regularity for f  is only dependent of translation (or rotation). When there is dependence of a trace in (u,v)
on a disk domain, this contradicts HE.

Assume $f(\zeta) \in H$ has negative exponential type. Consider continuation by a circular element, $f(\zeta + e^{v})$, where v is assumed independent of $\zeta$. We can using the Laplace transform $\mathcal{L}$, determine g so that $f(\zeta + e^{v})=\mathcal{L}(e^{-e^{v} s} g)(\zeta)$. Consider $<f,\phi>=< g,\mathcal{L}(\phi)>$, with $\phi \in \mathcal{S}$. Then we have $< e^{-e^{v} s} g, \mathcal{L} \phi>=<g, e^{e^{v} s} \mathcal{L} \phi>=<g, \mathcal{L}(\tau_{e^{v}} \phi)>=< \mathcal{L} g, \phi>(\zeta + e^{v})$. Assume $v \sim 1/\rho_{n}^{2}$, with $\rho_{n}$ real. Then we have given $\log \log g < 1/r_{n}$ and $\frac{1}{\rho_{n}^{2}} + \frac{1}{r_{n}} < 0$, that
$f(\zeta + e^{v})$ defines a continuation of negative type.
Note that given an analytic scalar product, we have monogenity according to reciprocal sub modules.
In particular $g(e^{e^{v}} \phi)=f(e^{v} + \log \phi)=f(e^{v} \log \psi)=h(v + \log \log \psi)$,
where $h=\mathcal{L}f$ and $\psi \in \mathcal{S}$.

We assume $\Omega$ a rectifiable surface, limited by a rectifiable boundary. Assume first derivatives to F analytic, are given by (X,Y) and $W=Y_{x} - X_{y}$, then

\newtheorem{lemma9}[lemma10]{Lemma}
\begin{lemma9}
For rectifiable sub level surfaces to F, with rectifiable first surfaces, we have $W \neq 0$.
\end{lemma9}

Let $f=F(\gamma)$, with F absolute continuous homeomorphism and $\gamma \in S^{1}$.
W=0 implies F real. Assume $\overline{F}(e^{\phi})=F(e^{\overline{\phi}})$. When $e^{\phi}$ on a first surface, the same holds for $e^{\overline{\phi}}$, in the same manner for the line between $\phi$ and $\overline{\phi}$, note that $dF \equiv 0$ on the line.
When $\{(x,y) \quad F=const \}$ is rectifiable and homogeneous, that is when y=y(x), we have that dy/dx $\neq const$, thus $\frac{d \overline \phi}{dt} / \frac{d \phi}{d t} \neq const$, that is not a line.

A SH function determines a rectifiable surface
(\cite{RieszF_26}). Consider u SH on rings, assume U the best harmonic majorant on rings =0 on the boundary, then we have, $I_{ij}=\int \frac{\delta U_{ij}}{\delta n} ds \geq 0$, where n the outer unit normal to curves in the ring. This is corresponding to $I_{ij} \sim \frac{I(r_{i})-I(r_{j})}{\log r_{i} - \log r_{j}}$.
$I(r)$  is a convex function of log r. When f is convex, $\{ (x,\lambda) \quad f < \lambda \}$ is a convex set. Convex sets in the plane define rectifiable surfaces. FU=f, with F SH and U harmonic, define a rectifiable surface.

ex : dU/dV=$\rho \rightarrow 1$ in $\infty$ regularily,
implies $\rho \notin L^{1}$, that is we may have spirals in limes, that do not occur on compact strata. Starting from Stieltjes integral,
that is $\Phi(U) f= \int f d U + \int_{\Gamma} df$, where $\int_{\Gamma(U)} df=0$ gives bi-linear action,
with interpolations property. When $0 \neq \int_{\Gamma} df$ linearity for $\Phi(U-V)f$, is dependent on
$\Gamma(U) - \Gamma(V)$.

ex: When (I) are rings of constant modulus, (I) can be given as kernel to a homeomorphism, $(I)=$ ker h. (x,y) defines a homogeneous space, when y=rot(x). Constant modulus for (x,y), is seen as the  boundary to a poly cylinder.

$\{ d \mid f \mid=\mid f \mid=0 \}$ is seen as a principally defined set. When sng is dependent on $(f,\overline{f})$, we have sets of higher order. Assume $(I)=ker h$, with $h=0$ implies $\mid f \mid=c$ and $dh=0$ implies f=c. Consider $h(f)=<H,f>=0$,
that is a set for which locally $H \bot f$. H preserves regularity, when H algebraic. ex : $<dH,f>=0$ implies f=c.

When F is regular in a homogeneous space (x,y), with y=y(x) regular, then F is regular in x. Reversible homogeneity implies that x can be determined from x=x(y).
ex: extremal ray according to $S(dU)=\{ dV \quad dU+dV=dI \} \subset G$, when $dI,dU \in G$.
On the other hand, outside an extremal ray, we may have $dI - dU \in S$ does not imply $d I \in G$. When S defines a solution to an undetermined moment problem,
we get infinitely many solutions. A determined moment problem, implies a unique solution.

\section{Subordinate elements}
Define a circular element, according to $e : \Sigma(z_{0},r)=\Sigma c_{n}(z-z_{0})^{\pm n}$, with $z \in D(z_{0},r)=\{ 1/\mid z-z_{0} \mid < r \}$. Assume $e' : \Sigma(z_{0},\rho)$, with $\rho=\min(\mid z-z_{0} \mid, r-\mid z-z_{0} \mid)$, then $\rho < r$ implies $e' \subset e$.
Existence of $e' \subset e$ and $e' \subset f$ implies f is continuation of e (\cite{Markush67}).

In particular a continuous group according to $e'$, can be continued to a continuous group according to e. ex : $dU \rightarrow dI \rightarrow dV$ regularly, defines a continuation of dU. Subordinate elements $e' \subset e$, correspond to e being divisible, when e irreducible
(order 0) it corresponds to presence of a boundary for convergence.
When E is a path, $\cup E_{i}$ gives a continuation of path. Relative the disk closure, E is seen as a subordinate element.

A d-dimensional regular set E, is called representable, if E is closed and if we have existence of $\mu$ Borel, with supp $\mu$=E and
$C_{0}^{-1} r^{d} \leq \mu(E \cap B(x,r)) \leq C_{0} r^{d}$.
Consider a very regular  curve C in 3-space and $\pi(C,p)$ the orthogonal projection on the tangent space. The geodesic curvature for C,  is the total curvature for $\pi(C,p)$, that is $K_{g}=\int_{0}^{r} k(s)ds$. Thus,
if $C=\Gamma_{1} \cup \Gamma_{2}$, we have that the curvature can be given by $k_{1} + k_{2}$, where we can have for instance $k_{1}=0$. Note that when we write $k=g'$, with g absolute continuous, the total curvature is given by $g(r)$. When g=0 for a curve on a rectifiable surface, $K_{g}(C)=0$ and further, when g is finite on a curve with finitely many segments, we must have $K_{g}(C) < \infty$.
ex : representable implies invertible, that is $dI \in G$ continuous and abelian. When $dI \notin G_{H}$, but $dU \in (G_{2})_{H}$, then $d U^{-1} \in (G_{8})_{H}$, that is $dU^{-1} \in (G_{2})_{C}$.

Given $V^{2} \rightarrow G_{2}$, with UV=VU, that is a disk B, the sub groups are seen as representable or not. Given representability, exponential type is preserved.
When $E \cap B$ is represented by a sub group, the set can be 0,1 or 2-dimensional, dependent of type. When the representation of sub group is not dense in B, we discuss codimension. Representabilty according to Radon-Nikodym, assumes a strict condition outside a compact set.

Algebraic continuation preserves geometric properties. dim E=dim $E \cap B$ is corresponding to a symmetry condition. Algebraicity is dependent of regularity for the continuation, for instance dUI=dIU, with an approximation property. Assume $dU \rightarrow d U^{\bot} \rightarrow dU^{\bot \bot} \hookrightarrow dU$, with $dU^{\bot}=dI-dU + dW$. Given $dW^{\bot}=dW$, this is not seen as a regular approximation and thus does not necessarily affect the dimension. A sufficient condition for absence of geodesics in the polar, is assuming the polar given in representation space, that we have rectifiability for curves in the polar. Given the polar regular (point wise topology),
G has algebraic action, corresponding to algebraic continuation in G. Given $f \in (I)_{H}(\Omega)$, we have for projective movements, that $dU(f) \in (I)_{H}(\Omega)$. Under a strict condition, we have density in $(I)_{L^{1}}(\Omega)$. We can determine a maximal sub group for algebraic action, in the sense of $C(dI)=\{ dU \in G \quad dUI=dIU \}$, relative exponential topology.

Assume f monogenous in $\xi$ and $\xi \rightarrow (u,v)$ absolute continuous. When f is monogenous, constant surfaces are negligible.
Over $S=\{ f=const \}$, $D_{S}(f)$ is finite, but $\tau_{S} \neq \tau_{\overline{S}}$, that is geodesic  curvature for S,  is not geodesic curvature for the stratification of the surface.
ex : $D_{1}(f)-T(f)=\int \sigma$, that is $\mid f' \mid^{2} - \mid f' \mid=\sigma$ as $\mid f' \mid < 1$, we have
$\mid f' \mid (1-\mid f' \mid) \sim \sigma$. When $\mid f' \mid=0$ or =1, we have $\sigma \sim 0$.
When $\mid f' \mid \rightarrow 0$ in $\infty$, we have $(1-\mid f' \mid) < \infty$ implies $\sigma \rightarrow 0$ in $\infty$, in particular integrability for $\sigma$ is dependent of the polar.

A maximal chain is defined by, non-existence of dV, with $dU \prec dV \prec dW$.
Regularly continuable implies $k_{g} \neq 0$ (non-trivial).

ex: Assume $x^{'2} + y^{'2} + z^{'2}=1$ (rectifiable), with $x_{j}' \in AB$ corresponding to BV. Then, presence of rectifiable trajectories, motivates existence of $x_{j}^{*'}$, according to $< x', x^{*'}>=1$. Cf removable sets, that is $\Sigma x_{j}^{'2} \in AB$ on E, implies $ \Sigma x_{j}^{'2}=const.$ on E, defines $\mathcal{O}_{AB}$ (\cite{AhlforsSario60}).

\section{Very regular boundary}
Assume A a perfect trace and B has isolated singularities, then $A \rightarrow_{ac} B$ gives continuation of a perfect set, that is (A,B) defines a very regular boundary. Regularity refers to the region outside (A,B).

Absolute continuity is not a radical property, we can consider $\Gamma_{0} \rightarrow_{\Phi} \Gamma_{1} \rightarrow_{\Phi} \Gamma_{2}$, with $\Phi^{2}$ absolute continuous. ex : $\Gamma_{1} \rightarrow \ldots \rightarrow \Gamma_{N}$ absolute continuous, defines a very regular boundary, when we have existence of $\Gamma_{N}$, with isolated singularities. ex:  bounded singularities with respect to dU.
Assume dI has point support and $(d U^{N}=dI)f=const$ implies $f=const$, in particular $df \neq 0$
implies $d U^{N}f \neq dI f$ or $\frac{d^{m} \Phi}{d t^{m}} f \neq 0$, corresponding to meromorphic representation. When $d U^{N}=dI$ defines a rectifiable boundary, according to $\Phi^{N}$ absolute continuous, this does not imply that $dU=dI$ is rectifiable, given dU BV.

ex: Consider $\Gamma_{1} \cup \Gamma_{2}$, with $\Gamma_{2}$ rectifiable and
$\Gamma_{1} \rightarrow \Gamma_{2}$ absolute continuous.
Given a very regular boundary according to $\int_{\Gamma_{1} \cup \Gamma_{2}}=\int \rho_{1}\rho_{2} dx$,
we consider irreducibles according to $\rho_{1}=0$ implies $\rho_{2} \neq 0$, that is the multipliers do not have common zeros. Further $\rho_{1}=const$ implies $\rho_{2}' \neq 0$.
Note that a projective transformation does not imply preservation of singularities. Ex: $dU=\rho dI$, with $\rho=0$ non-trivial.

ex : When $(dU,dV)=(\rho dI,\sigma dI)$, with $\rho=const$ orthogonal to $\sigma=const$, we can, by deformation, assume $\rho=\rho(x)$ and $\sigma=\sigma(y)$. When $\rho,\sigma$ holomorphic approximations, we can use the inverse representation (\cite{Oka60}), assume $S=\{ F(x_{1},\ldots,x_{n};\rho_{1},\ldots,\rho_{m})=0 \}$ a homogeneous space, with $\rho_{j}(x)$ holomorphic and such that
$(x,\rho) \rightarrow x$ 1-1. Assume $\sigma$ a ramification surface (order $\mu-1$), with irreducible component $\sigma_{0}$. By a pseudo conformal transformation $\sigma_{0}=\{ x_{1}=0 \}$. Then $\rho_{1} =a_{0} + a_{1} t + \ldots$, with $t=x_{1}^{1/\mu}$.

\section{Principal operators}
Assume $P= P_{m} +R$, with $R \prec \prec P_{m}$.
Starting from $C^{-1} \mid \xi \mid^{a} \leq \mid P_{m}(\xi) \mid \leq C \mid \xi \mid^{a}$, then $P_{m}=0$ implies that $R=0$ (outside a compact). In particular $P_{m}=const$ implies $\mid \xi \mid=const$, that is a rectifiable surfaces outside compact.
Conversely, we can prove, when $P_{m}=c$ is rectifiable, $\{ P=c \} \rightarrow \{ P_{m}=c \}$ absolute continuous. In particular $m \{ P=c \}=0$ implies $m \{ P_{m}=c \}=0$. ex : Assume $P=\alpha P_{m}$, where $\alpha$ is locally reduced (downward bounded). Then $\{ P < \lambda \} \rightarrow \{ P_{m} < \lambda \}$ is continuous. Given P HE and $\alpha$ bounded, we have $P'/P \sim \alpha' P_{m}/P \rightarrow 0$. When P (not HE) is polynomial, we only have $\alpha'$ bounded. Further, given $\alpha'$ single valued $\alpha$ is monogenous.

For f with $f^{(p)} \neq 0$,  we have existence of V(f) polar, where f is regular outside V(f). Consider continuation $T_{\gamma}(dU)$, with for instance $dU=\rho dV$ and $dV$ rotation. Then $T_{\gamma}^{k}f$ is regular, with bounded singularities relative dU, when $dU^{k}$ regular on supp f. Assume further a projective resolution and $(T_{\gamma}^{k} - Id)^{m} \sim T_{\gamma}^{m}(dV_{k})$ regular, this implies singularities bounded relative $dV_{k}^{m}$.

Assume U harmonic and $f$ SH.
Then the composition f(U) defines an absolute continuous mapping in U (\cite{Montel28}). In particular, if $(\Gamma_{m})=\{ P_{m}=const \} \subset \{ U=const \}$ and f a SH mapping between boundaries, when f=const on $\Gamma_{m}$, it follows that f is absolute continuous in U.
$\frac{df}{dU}=0$ implies $f(U)=const$. $ \frac{d}{dt} U=0$ as $U=const$, that is $\frac{df}{dt}=0$ and f absolute continuous preserves constant surfaces. Convex sets are path wise connected, with an interpolations property.

\newtheorem{lemma6}[lemma10]{Lemma}
\begin{lemma6}
 $\{ P=c \} \rightarrow \{ P_{m}=c \}$ is absolute continuous, given P HE.
\end{lemma6}

Given P HE, $P_{m}$ must be HE, thus $\{ P_{m}=c \}$ are rectifiable according to an earlier lemma. The boundary homeomorphism (\cite{Collingwood66}) identifies 1-1, the first surfaces above.
$P_{m}(tx)=t^{m}P_{m}(x)$ implies given $S=\{ P_{m}(x)=c \}$ and $x \in S$, when $t x \in S$, that t must  be finite, that is $P_{m}=c$ normal and rectifiable.

Assume $\Sigma F_{j}^{2}$ defines an orthogonal (projective) decomposition.
Triangulation is dependent of the polar. Assume $S=\{ P=c \}$, such that S has regular holomorphic neighborhoods locally.
Assume P=$P_{m} + R$, with $R \prec \prec P_{m}$. Thus, $P/P_{m}=1 + R/P_{m} \rightarrow 1$ in $\infty$. P=const implies $P_{m}=const$.

To determine spiral geodesics, it is sufficient to consider a second order boundary condition, that is m=2.
h harmonic, corresponds to a real transformation. $P_{m}$ defines a sub module.
ex : $< f,dU>$ absolute continuous, implies Uf=const iff $d < f,dU>=0$, further
$< f(d \phi),dU>=<f,dU>(d \phi)=0$, that is $Uf(d \phi)=0$.

\section{Hypo invariants}
(I) invariant for G, if $G(I) \subset (I)$. That is, $\forall dU \in G$, $dU(f)=dI(f)$
and $f \in L$ defines L as an invariant ideal. Assume $d G_{3}=G_{2}$, has a continuous inverse, then a line in $G_{2}$ corresponds to line in (x,y). L fixed defines a boundary in (x,y,z) corresponding to a hyperplane.
Starting from $f(u) \rightarrow f(0)$ regularly and regular transformations between linearly independent movements in a continuous group G, we get an approximation property for G.

When $(I)(\Omega)$ is holomorphically convex, we do not have existence of dV, a continuation of R(G), that is we do not have perfect transformations, relative (I).
Note that the concept is relative topology, that is $C(G_{2}) \simeq H(G_{8})$, for instance $G_{2}$ is not  continuable in H, but continuable in C. Given a regular contraction, we have that $\gamma \in R(G_{8})$ has a representation in $R(G_{2})$. Perfect transformations have codimension, non-continuable groups
correspond to maximal chains.

D(f) induces a semi-norm on $\mathcal{D}(\Omega)$.
Let $D(f,g)=\int f' \overline{g}'$. ex : $g'=\widehat{f}'$.
$\phi$  is hypo continuous (\cite{Treves67}) if, every nbhd 0 in G, has a nbhd 0 in E,F, such that $(x,y) \in$ nbhd (0,0) implies $\phi(x,y) \in$ nbhd 0. In particular $\tau(\phi(x,y)) \leq p(x) q(y)$, where p,q,$\tau$ are semi norms on E,F and G. Further,
that $\phi_{x}(y)$ and $\phi_{y}(x)$ are equicontinuous, for x,y bounded.
ex : $ \frac{d}{dt} (\gamma,f(\gamma))$, that is $(\frac{d \gamma}{d t}, \frac{d f(\gamma)}{d t})$ single valued, means in particular $df/d \gamma=0$ implies $d \gamma / dt=0$ (negligible).

 ex: xq-yp=0 implies q/p=y/x. That is $1/k < \mid dy/dx \mid < k$ implies $1/k < \mid y/x \mid < k$. If further xp+yq=0,
we have that $1/l < \mid x/y \mid < l$. ex : $x \rightarrow y$ linear.
$1/ k  < \mid dy/dx \mid < k$ implies $1/l < \mid dx / dy \mid < l$, that is a reversible condition that does not exclude traces.

Given $dI \in G$, with $dUf=dIf $, $\forall dU \in G$ or $f(u_{1},\ldots,u_{8})=f(0,\ldots,0)=const$, that is a higher order diagonal (zeros). Consider $dU^{-1}V=dI$, that defines the diagonal relative conjugation, for instance using $(dU,dV) \rightarrow (dU,dV^{-1}U)$. For a very regular group, $dV^{-1}U \rightarrow d U^{-1}V$ continuous and $\lim V^{-1}U=\lim U^{-1}V$. ex : f=const iff $(u,v) \in \Delta$, where $\Delta$ defines trace, defines a very regular domain.

When $\phi(x,y) =x-y$, we have that a strict condition implies density for translates, $\phi^{-1}(0)=$ nbhd (0,0), that is presence of trace implies $\phi$ is not hypo-continuous  mapping. The same holds for its convex closure.
Further, $(x,y) \rightarrow (x,y/x)$ is not single valued, when $y/x=c$.  Presence of invariants sets (lineality) gives a non-monogenous continuation, that is multivalentness for limes. Note that for $\Delta$ translation invariant sets, we have $\Delta(df) \neq \Delta(f)$. When $\Delta(f)$ non-trivial, we do not have a single valued representation. Given a regular contraction, $G \rightarrow G_{2}$, we have that monogenity implies absence of translation invariant sets.
Further, when $d V^{-1}U=dI$ defines the trace, we have that $d V^{-1}U \in G$ implies G is discontinuous.

A convex set in the plane,  is one sided relative a plane.
$D(f) < \infty$ does not imply a convex set or $f \rightarrow 0$.
ex : $\frac{\delta^{2} f}{\delta x \delta y} f(x) g(y)=f'(x) g'(y)$. Note that $< dH,f' g'> \sim < H, f'' g''>$, given $dH=dI(x) \otimes dI(y)$

A rectifiable surface, with coordinate functions of separate BV, is hypocontinuous. A surface of finite D-semi-norm may have trace sets in infinity, thus it is not necessarily hypocontinuous.

\section{Geodesic curvature}
In the representation space, we assume a contraction property $G \rightarrow G_{2}$. For a monogenous function, we have $dN \bot G$ implies $dN \bot G_{2}$. For a rectifiable surface, the curvature is determined as the total  curvature for the orthogonal projection on the tangent space (\cite{Backlund16}). Thus, for a rectifiable surface, geodesics can be determined in the representation space. The geodesic  curvature is given in the plane by 1/r, where r  is the radius to the osculating circle. dU=dV implies 1/r=0.
Note that monogenous on representation space implies decomposable. When $V^{n-1}$ is a cylinder web, we have spiral geodesics. Given a local isometry according to Brelot and when (dU,dV) represents a very regular domain, geodesic  curvature is preserved in (x,y,z).  A rectifiable surface can include curves with geodesic  curvature =0, but the total curvature for rectifiable curves of positive curvature, is not =0.

For a surface M in 3-space, along a curve $\gamma$ on M, consider the tangent T to $\gamma$, the normal N to M along $\gamma$ and form T $\times$ N (cross-product). Assume n normal curvature to M and g geodesic curvature, then we have $\frac{d T}{ds}=g (T \times N) + n N$.
ex : dU+dV=dI defines dN locally. Given a spiral $\gamma$ on a cylinder web in 3-space, we see that $\gamma''$ describes the normal to S, why it has geodesic curvature =0.

Bäcklund gives a deformation, where two surfaces can be deformed on each other, if zero curves are mapped on zero curves (projective  mapping) and if under the deformation,
we have preservation of length and geodesic curvature for the curves (\cite{Backlund16}).
When the surface is rectifiable in the sense that the surface can be generated by its tangents, Bäcklund proves by deformation, that geodesic curvature for the surface can be determined by geodesic curvature for envelop of tangent planes to the surface.

The deformation of integration path (\cite{Nils80}) starts from algebraic $\xi_{j}(y)$, in $V=\{ p(y,x)=0 \}$, with y close to $\infty$, and their Puiseux developments. Rings $\Omega(\alpha,K)$, with modulus K, are formed by identifying leading coefficients in the development, relative modulus and exponents. Points in rings satisfy $(K-\epsilon) \mid y \mid^{\alpha} \leq \mid x \mid \leq (K+\epsilon) \mid y \mid^{\alpha}$, as y close to $\infty$. $\Omega(\alpha,K)$ are disjoint for y close to $\infty$.
Groups $\Omega_{s}$ are defined by point wise equal leading coefficients.
Every group can be divided into second order rings, starting from second leading coefficient, and so on.
In this manner all $\xi_{i}$ are situated in different rings (groups).
The rings are combined by straight line segment, to a connected path.

Equivalence relative a scalar product, according to $<Uf,Vg>=<f,g>$, assumes a rigid geometry. When we include a reference to flux according to
$<Uf,Vg>-<f,g>=\int_{\Gamma(u,v)}$, we assume equivalence is independent of the boundary condition. $V \sim U$ implies $\int_{\Gamma(u,v)}=0$. Dependence of the boundary gives an elastic model.
ex : $\{ f^{j}=c \}=\Gamma_{j}$, where $\Gamma_{N}$ is a rectifiable zero line. Assume $\Gamma$ is mapped on $V^{n-1}$, then we have $\Gamma \rightarrow \Gamma_{N}$ is absolute continuous and $(\Gamma,\Gamma_{N})$ defines the boundary.

Vanishing flux is defined by $\int_{\Gamma} df^{\diamondsuit}=0$, with $df^{\diamondsuit}=\alpha dx + \beta dy=0$ on $\Gamma$ iff dy/dx=$-\alpha/\beta$. df=0 on $\Gamma^{\diamondsuit}$ iff $dy/dx=\beta/\alpha$. $\int_{\Gamma} df$ is interpreted as $\Gamma$ approximates $\{ f=c \}$ according to $z \rightarrow -1/z$.

Assume $D(f) < \infty$, with f polynomial on strata. Consider an absolute continuous (non-constant) deformation $\tilde{S}_{j}$ as above, of strata $S_{j}$.
Thus $f'(x,y)$ single valued and finite (monogenous) on $\tilde{S}_{j}$. Note that we do not have that
geodesic  curvature is preserved on $\cup \tilde{S}_{j}$, unless $\cup S_{j}$ is rectifiable as surface and consequently hypocontinuous.

ex : Consider $f + \lambda g$ ,with g a linear form, for $\lambda$ scalar, invariant for monodromy. Assume in particular $f + \lambda g \simeq f(\phi(x,y))$, where $\phi$ a transformation over $G_{2}$. Starting from critical points, the deformation multivalent.
When $\phi$ is BV, $f + \lambda g=c$ is rectifiable .

Monogenity over G preserves continuity for G, that is $dI(u,v) \in G$ implies $dI(u) \in G$. Sufficient for this is $\Gamma(u,v) \rightarrow \Gamma(u)$ continuous, with $\Gamma$ rectifiable. Necessary is $\Gamma_{1}(u,v) \rightarrow V^{n-1}$ according to Brelot, separates points, in particular maps on to a 2-boundary. $\Gamma_{1}(u,v) \rightarrow V^{2}$ with (u,v) very regular, preserves monogenity, given $V^{2}$ rectifiable. Further, sets of measure zero are preserved.

Monogenity for F(u,v) is sufficient for ``transversal'' single valued normal. A very regular boundary $\{ \Gamma_{j} \}$, defines multivalent normal. $\Gamma_{1} \rightarrow \Gamma_{2}$ absolute continuous, implies a connected domain.

ex : $f=g'$ with $\int_{E_{n}} g'(y,x)dx=\int g(y,x) dI(x)$. We consider a deformation of $E_{n}$, by a very regular group. $<Vg,dU>(\phi)=<g(\phi),dV \times dU>$.  The implication $g \in H(u)$, is dependent of the boundary condition, $\int_{\Gamma(u,v)} dg=0$, that is independence of V. $\Gamma \downarrow \{ 0 \}$ relative continuation, is interpreted as bounded singularities relative (dU,dV).

ex : $C=\{ u=v \}$ $C'=\{ \mid u-v \mid=R \}$. Assume geodesic curvature =0 on C. An absolute continuous  mapping to $C'$, form a very regular boundary with positive geodesic curvature. Assume $p,q \in$ nbhd C and p is mapped on a curve with positive geodesic curvature in nbhd $V^{n-1}$, then also the image of q has positive geodesic curvature. Local isometries (Brelot) preserve geodesics, however the boundary mapping is not necessarily 1-1.
Note that $C \rightarrow C'$ absolute continuous, it is not assumed a reversible property.

ex : $dU+dV=dI - dW$. Reflexivity implies $dW=d W^{\bot}$, which defines the polar. $dW=0$ implies $dV=d U^{\bot}$.
When $R(dU,dV)^{\bot}=R(dW)$. $dW \bot dI$ implies absence of trace.
dU=-dV implies dW=dI.

\newtheorem{ex2}[lemma10]{Example}
\begin{ex2}
Assume $\Gamma=\{ z \quad f(-1/z)=const \}$ rectifiable in the plane and assume $S^{1}$ an osculating circle, $\tau_{g} \Gamma=\tau_{g} S^{1}=1/r$. Then we have that $\int_{\Gamma} df=\int_{\Gamma} \frac{\delta f}{\delta n} \mid -d 1/z \mid$ (\cite{AhlforsSario60})
\end{ex2}

Starting from a curve of unit speed and a standard frame in 3-space: T,N,X=$T \times N$, we can give g= geodesic curvature, k=normal  curvature and
t=principal  curvature. Assume X/N=$\mu$, T/N=$\rho$, where $\mu,\rho$ are assumed real, then we have $T'/T=g \mu/\rho + k/\rho$ and $X'/X=-g \rho/ \mu + t/\mu$.  ex : $d \log T=d \log X$ and k=t=0, then g=0. When k=t=0, the curve is flat along with the tangent plane.

Starting from $x(u,v)$ monogenic, when $du/dn=\rho$ and $dv/dn=\mu$, we assume $dv/du =\mu/\rho \in L^{1}$ implies $\mu,1/\rho \in L^{1}$. When dv/dn=const, it is in $L^{1}$ on compact sets. With these conditions, monogenity locally approximates a frame where the multipliers are in $L^{1}$.

When we have absence of spiral traces, f(u,v) can be considered as f(u+iv). Cauchy Riemann gives f(u+iv)=$f_{1}+i f_{2}$ and
$f'(u+iv)=d f_{1}/d u + i d f_{2}/ du$,
that is it is sufficient to consider rotation.
The pre image to the Brelot homeomorphism, is homogeneous v=v(u). Monogenity according to $<u,<v,n>>=<<u,v>,n>=<v,<u,n>>$, means that orthogonal projection preserves homogeneity for the space.

\section{Analytic stratification}
Given a compact manifold M, for a sub module L of vector fields over analytic functions, define integral elements to L in x, $L(x)=\{ u(x) \quad u \in L \} \subset T_{x} M$ and the isotropy sub algebra $L_{x}=\{ u \in L \quad u(x)=0 \}$, then
$0 \rightarrow L_{x} \rightarrow L \rightarrow L(x) \rightarrow 0$ is exact (\cite{Nagano66}).
In particular, when dI evaluation in x and $dI,dU \in G$, we have by exactness $dUI=dIU \in G$.
Theorem: $M=\cup N$ uniquely, with $N \subset L$ (\cite{Nagano66}).

Assume $d L_{k}$ define convex domains according to $(d L_{k})$ inner/outer relative $d L_{k}$. $d L_{k}$, is assumed for a stratification to be rectifiable. In particular, $d L_{k}=\Phi(S^{1})$, with $\Phi$ absolute continuous. We assume the stratification approximates first surfaces S, by locally convex sets. Given multivalent boundary (constant dimension) we assume an absolute continuous boundary mapping,
in particular, $\int_{\Phi(S^{1})} df=\int_{S^{1}} d f \circ \Phi^{-1}$.

Assume $bd L_{j} \supset L_{k}$, for some $k \neq j$. Further, assume a Brelot homeomorphism h, that preserves the stratification,
with $V^{n-1}=h(bd G_{n})$ rectifiable. We assume as a condition for stratification, the normal is single valued. A planar tangent surface in 3-space, defines a single valued normal.
ex : $V^{-1}U \bot W$,$U^{-1}V \bot W$ and UV=VU. When V=U, we have $I \bot W$, that is $I/W \rightarrow 0$ in $\infty$.
For instance W reduced according to $W \phi=0$ implies $\phi=0$. Note that U=V is not seen as rectifiable.

ex : $G_{2}'=(dU,dV) \rightarrow (dU,d V^{-1}U)$ (ähnlich). $dV^{-1}U=dI$ identifies an extremal evaluation in G, only when G is continuous.
$dV^{-1}UV \sim dU$ and $d U^{-1}V U \sim dV$ implies $ Uf= \int Vf d V^{-1}U= \int f dU$.
An abelian group does not imply a continuous group. ex : $dV^{-1}U \in G$ but $dI \notin G$. $dV,dU \in G$ with $dI \notin G$ does not imply $dV^{-1}U \in G$.

When f is monogenic, we can assume D(f)=0 implies T(f)=0. Since zero lines are mapped on zero lines, the corresponding mapping is projective. Further, an monogenic continuation $\{ T(f) < \infty \} \rightarrow \{ D_{1}(f) < \infty \}$  is projective. In general, $T(f) \rightarrow D(f)$ defines a discontinuous continuation, for instance $\tilde{f}$ BV does not imply $\tilde{f}$ analytic (or continuous). ex : f analytic and $D_{1}(f)(G_{2}) < \infty$ implies $T(f)(G_{8}) < \infty$, that is without approximation property in $G_{2}$. An analytic stratification in $G_{2}$ has a contraction property.

ex (Riesz-Thorin) : $L^{1} \cap L^{2} \subset L^{\alpha(\theta)} \subset L^{1} + L^{2}$,
where $1/\alpha(\theta)=1-\theta/2$. Thus, $T(f),D_{1}(f) < \infty$ implies $\int \mid f' \mid^{\alpha(\theta)} < \infty$.

Assume $M_{\lambda}=\{ f < \lambda \}$ is limited by a rectifiable surface $\{ f= \lambda \}$, for instance $dM_{\lambda}=f(S^{1})$ with f absolute continuous homeomorphism. Sufficient for an algebraic polar V(f), is strict pseudo convexity.

$e(\lambda)=\int_{M_{\lambda}} g(\xi) d \xi$, with for instance $M_{\lambda}=\{ f(\xi) < \lambda \}$
defines an orthogonal projection operator $E_{\lambda}$, according to $<E_{\lambda} \varphi,\psi>=<e(\lambda), \varphi \otimes \psi>$, for instance a spectral function related to  $M_{\lambda}$, where $g(x,y;\xi)=(2 \pi)^{-n} e^{i(x-y) \cdot \xi}$ (\cite{Nils72}). When f HE relative a continuous group $G_{H}$, we have monogenity according to $(A+B)^{\bot}=B^{\bot} A^{\bot}$, HL=$\{ 0 \}$ iff $A \bot B$.

\newtheorem{lemma19}[lemma10]{Lemma}
\begin{lemma19}
 The spectral projection operator defines a $C^{\infty}$ stratification, given the symbol HE.
\end{lemma19}

Assume $f \in C^{\infty}$ (or real-analytic), satisfies a strict condition in for instance $L^{1}$.
Let $M_{\lambda}=\{ \lambda_{0} < f(x) < \lambda \}$ be compact in $\mathbf{R}^{n}$, for every real $\lambda$. We consider a spectral kernel $e(\lambda)=\int_{M_{\lambda}} g(x) dx$, for a regular function g, that defines a spectral projection operator.

$S_{\lambda}=\{ f(x)=\lambda \}$  is a $C^{\infty}$ manifold, (n-1)-dimensional and defines a (n-1)-cycle (\cite{Nils80}). $d e(\lambda)/d \lambda=\frac{d}{d \lambda} \int_{M_{\lambda}} g(x) dx$, as $d \lambda=d f$ on $S_{\lambda}$,
by Stokes, $d e(\lambda)/d \lambda=\int_{S_{\lambda}} g dx / d f$.

In particular, as $\varphi \sim < x-y,\xi>$, we can consider $e(\lambda)$, as an asymptotic integral. Further, (\cite{Nils72}), the spectral function to a formally HE differential operator, with variable coefficients,  is regularizing in x,y and satisfies
$C^{-1} \lambda^{a} (\log \lambda)^{t} \leq e^{(2 \alpha)}(\lambda) \leq C \lambda^{a} (\log \lambda)^{t}$,
when $\lambda$ large, $C>0$, t entire number and $a \geq 0$.

Assume $f=\Sigma g_{j} \frac{\delta \varphi}{\delta x_{j}}$ (\cite{Malgrange73}), with regular coefficients, flat on $S(\varphi)$, that is f is defined relative involution. Thus, $\varphi$ determines the type of movement. Movements are assumed as linearly independent over strata in the representation space. Monodromy determines sub groups for regular continuation.

Assume X reduced, is defined by $f_{j}=0$ for $j=1,\ldots,k$ and $f_{j} \in H$. ex : $f=P/Q$ for P,Q polynomial, with Q reduced, that is $Q \in X$ and
Qf polynomial. A very regular boundary can be defined by $X'=\{ \exists f_{j}=0 \}$. Further, $fg=0$ defines $X'$;  $\mid f \mid + \mid g \mid=0$ defines X.

\newtheorem{ex3}[lemma10]{Example}
\begin{ex3} Assume $f,g \in C_{0}^{\infty}$, then
H(f+g)=0 does not imply H(f) + H(g)=0, even when $\delta(f+g)=0$ implies $\delta (f) + \delta (g)=0$.
When a geometric ideal (I) is generated by $f,g$ with $dg/df \neq const$, such that $f,g$ have support on $\Omega$
implies f+g has support on $\Omega$, then H gives an extremal definition of $(I)(\Omega)$.
\end{ex3}

Assume $\gamma_{j}$ regular and $f_{j}=<\gamma_{j},d \mu>$, with $d \mu$ a reduced measure and $d \mu \rightarrow dI$.
Then $\delta(\gamma_{1} - \gamma_{2})=0$ implies $\gamma_{1}(0)=\gamma_{2}(0)$. If
$<H, d (f_{1}-f_{2})>= <\delta, f_{1} -f_{2}>=0$, we have $< \gamma_{1}-\gamma_{2},d \mu>=0$ implies $\gamma_{1}=\gamma_{2}$. Given $d \mu$ an analytic measure (projective for harmonic conjugation), we have that $\gamma_{1}=\gamma_{2}$ implies $\gamma_{1}(0)=\gamma_{2}(0)$.
H is a separating functional, but $\delta$ is not.

Consider $Uf=\int f dU + \int_{\Gamma} df$. That is $Uf - f=\int f (dU-dI) + \int_{\Gamma-0} df$. $Uf \rightarrow f(u)$, with u single valued, assumes monodromy.
Assume $\Gamma$ multivalent and $\Gamma_{j}$ first surfaces $\not\ni \infty$, with $\Gamma_{1} \rightarrow \Gamma_{2}$ absolute continuous. Then we have df=0 on $\Gamma_{2}$ implies df=0 on $\Gamma_{1}$, that is a projective mapping ($ \in G$).
The representation can be extended to $f(\phi) \in C^{\infty}$, analogous with vector valued distributions,
 that is $<f,dU>(\phi)$, where dU is considered in $\mathcal{E}^{'(0)}$.
ex: $\Gamma^{\diamondsuit} \subset S_{\lambda}=dM_{\lambda}$. Given vanishing flux, we have $\int_{\Gamma} df(-1/z)=0$. Given $\Gamma$ en entire line, we have $z \in \Gamma$ iff $-1/z \in \Gamma$.

Starting from $G_{8}$ and $dU_{1} + \ldots + d U_{8}=dI$ a local representation, relative an irreducible algebraic equation F=0 (polar definition),
there is a parametrization in dN according to Puiseux, that includes a spiral axis. Thus, spiral solutions are
necessary for density. If we have absence of monogenic solutions, there is presence of spiral solutions.
More precisely, assume $F(u,n)=p_{0}(u) + \ldots + p_{n}(u) n^{k} \equiv 0$,
where $p_{j}$ are polynomial, without common zeros and $p_{k} \not\equiv 0$ of degree k. Assume $u_{0} \notin$ sng F and starting from zeros to $p_{k}$, the solution $n=\eta$ can be determined as a Puiseux series $\eta(u)=\Sigma a_{i} (u-u_{0})^{i/k}$, with $j \leq k$ (\cite{Markush67}). Monogenic contraction is sufficient to
exclude spiral normals.

\section{Monodromy}
Monodromy defines $(\gamma,f(\gamma))$
as single valued, where $\gamma$ is seen as connected. Thus, for f analytic, $f(G_{2})=f(G_{2}')$ implies $G_{2} \simeq G_{2}'$. Define $(I)=\{$ paths between two points $\}$. Given $<f(\gamma),dU>=<f(\gamma),dV>$, for all $\gamma \in (I)$, we can identify dU=dV in the measure space.

Assume $f \in (I)$, with topology in $\mathcal{D}_{L^{1}}$. Then we have $Uf - f \in \mathcal{D}_{L^{1}}$ identifies movements in representation space. That is, Uf=$\int f dU + \int_{\Gamma} df$ and Uf-f=$\int_{\Gamma \backslash 0} df$. When G is continuous and $\int_{\Gamma \backslash 0} df=0$ in $\mathcal{D}_{L^{1}}$, further when dU is 1-parameter, we can determine $dU \in G$.

$\Phi(u,v)f=\int f du dv + \int_{\Gamma} df \rightarrow f$ implies $\int_{\Gamma(u,v)} df=0$,
when the representation space is reduced to a 2-sub group (u,v). HE implies monogenity over $G_{2}$, why we can identify $G_{2} \simeq G_{1} \times G_{1}$. This assumes $\int_{\Gamma(u)} \rightarrow \int_{\Gamma(v)}$ projective and vice versa.

Given a strict condition, for instance $(I_{1}) \subset (I_{2})$, projective continuations can be determined.
ex : Reciprocal sub modules with $(M+N)^{-1}=M^{-1}N^{-1}$, where an extremal ray defines a reciprocal (monogenous) sub module.
More precisely, $(M+N) \in G^{-1}$ implies $M,N \in G^{-1}$, implies $M N \in G^{-1}$. Assume C defines a continuation of dM, then
$C^{*} \supset \{ <H,dN>(f) = < \delta, N>(f) \geq 0 \quad \forall dN \in C \}$, that is $ H \in C^{*}$, the canonical dual to C, analogous with one sidedness.

\newtheorem{lemma2}[lemma10]{Lemma}
\begin{lemma2}
Monodromy for the spectral projector to f HE, identifies very regular sub groups.
\end{lemma2}

$\mathcal{D}_{L^{1}}(X \times Y) \simeq \mathcal{D}_{L^{1}}(X) \otimes \mathcal{D}_{L^{1}}(Y)$, when X,Y are open sets in for instance the representation space. Assume restriction to a 2-subgroup, with an approximation property, that is $Id \in \overline{\mathcal{L}_{c}}(G,G_{2})$, with uniform convergence on convex, compact and equilibrated subsets. Then the 2-subsets can be identified as above. Given conjugation and monogenity, we can further identify the conjugated 1-subsets. Note that when the problem is to determine type of movement, it is sufficient to consider $\log f \in \mathcal{D}_{L^{1}}$.

Consider $G \rightarrow G_{2}$, a contraction to a planar model. When $G_{2}$ is very regular, the regularity for the symbol f i preserved. Necessary for HE is further monogenity over $G_{2}$.  When the group is abelian, the condition is symmetric. ex : $f(u,v)=(\rho_{1}/\rho_{2}) f$, is such that
$\lim (\rho_{1}/\rho_{2})f=\lim (\rho_{1} \rho_{2}) f$, as $\rho_{2} \rightarrow 1$. Further, $\lim f(u,v)=\lim f(v,u)$.
ex : A sufficient condition is $M=\{ \rho,1/\rho \in L^{1} \}(\Omega)$ with $\Omega$ unbounded and dU/dV=$\rho$.
Note that as f HE, we have absence of essential singularities in $\infty$. Over M we have an interpolation property for $G_{2}$.
For instance, continuation by N(u,v), under preservation of monogenity, that is N(u,0) preserves regularity.

Assume monogenity for f over $G_{2}$ and consider continuation $G_{2} \rightarrow G$ under preservation of monogenity,
sufficient for this is absolute continuity (non-constant), in particular deformation of a rectifiable surface, under preservation of length of curve segments and geodesic curvature (\cite{Backlund16}).
Starting from a projective decomposition (triangulation) dU+dV=dI, by decomposability $f(u,v)=f_{1}(u) \times f_{2}(v)$, we can identify $f_{j}$, given a very regular sub group. Assume $E_{\lambda}$ the spectral projection operator to a HE f. When $\widehat{E_{\lambda}}$ is corresponds to the characteristic function to sub level surfaces $M_{\lambda}$, we have that the kernel to E has regularization action, why (I-E) has very regular action and we have an approximation property. As long as (u,v) preserves compact sub level surfaces, the regularity for $E_{\lambda}$ is preserved.

Consider $\Gamma=\{ \Gamma_{i} \}$, a trajectory in a rectifiable surface. Then we can determine the
geodesic curvature, as $\tau_{g} \Gamma=\Sigma_{j} \tau_{g} (\Gamma_{j} \rightarrow G_{j} \rightarrow h(G_{j}))$, where $\Gamma_{j} \rightarrow G_{j}$ can be given by an orthogonal projection, for $\Gamma_{i}$
rectifiable  and h is the Brelot mapping. Presence of trace ($\tau_{g} \Gamma_{i}=0$) corresponds to a discontinous group, when $\Gamma$ rectifiable. Monodromy corresponds to a regular approximation property, relative G very regular, for instance $\Gamma_{1}$ rectifiable and $\Gamma_{i} \rightarrow \Gamma_{1}$ absolute continuous.

\bibliographystyle{amsplain}
\bibliography{januari_26}

\end{document}